\documentclass[10pt,a4paper]{article}
\usepackage{fullpage}
\usepackage[pdftex]{graphicx} 
\usepackage{amsmath, amsfonts, amssymb, amsbsy}
\usepackage{latexsym} \usepackage{stmaryrd}

\usepackage[labelfont=bf]{caption}
\usepackage{hyperref}
\hypersetup{
    colorlinks,%
    citecolor=black,%
    filecolor=black,%
    linkcolor=black,%
    urlcolor=black
}

\begin{document}

\newtheorem{theorem}{Theorem}[section]
\newtheorem{proposition}[theorem]{Proposition}
\newtheorem{corollary}[theorem]{Corollary}
\newtheorem{lemma}[theorem]{Lemma}
\newtheorem{definition}[theorem]{Definition}
\newtheorem{remark}[theorem]{Remark}
\newtheorem{example}[theorem]{Example}

\newcommand{\proof}{\noindent \textbf{Proof. }}
\newcommand{\qed}{ \hfill {\vrule width 6pt height 6pt depth 0pt} \medskip }
\newcommand{\separe}{\medskip \centerline{\tt -------------------------------------------- } \medskip}
\newcommand{\doit}[1]{\hfill {\small \tt [#1]}}
\newcommand{\note}[1]{{\footnotesize #1}}

\newcommand{\eps}{\varepsilon} \renewcommand{\det}{\mathrm{det}} \newcommand{\argmin}{ \mathrm{argmin} \,}
\newcommand{\Om}{\Omega} \def\interior{\mathaccent'27} 
\newcommand{\weakto}{ \rightharpoonup}  \newcommand{\weakstarto}{\stackrel{*}{\rightharpoonup}}
\newcommand{\R}{\mathbb{R}}

\newcommand{\stress}{\boldsymbol{\sigma}} \newcommand{\strain}{\boldsymbol{\epsilon}} 
\newcommand{\neu}{ \partial_{\mbox{\it \tiny N}} } \newcommand{\dir}{\partial_{\mbox{\it \tiny D\hspace{0.5pt}}} }
\newcommand{\jump}[1]{\llbracket #1 \rrbracket}
\def\Xint#1{\mathchoice
 {\XXint\displaystyle\textstyle{#1}}%
 {\XXint\textstyle\scriptstyle{#1}}%
 {\XXint\scriptstyle\scriptscriptstyle{#1}}%
 {\XXint\scriptscriptstyle\scriptscriptstyle{#1}}%
 \!\int}
\def\XXint#1#2#3{{\setbox0=\hbox{$#1{#2#3}{\int}$}
 \vcenter{\hbox{$#2#3$}}\kern-.5\wd0}}
 \def\ddashint{\Xint=} 
 \def\dashint{\Xint-}

\newcommand{\F}{\mathcal{F}}\newcommand{\wF}{\tilde{\mathcal{F}}}
\newcommand{\I}{\mathcal{I}}\newcommand{\J}{\mathcal{J}}
\renewcommand{\H}{\mathcal{H}}
\newcommand{\U}{\mathcal{U}}
\newcommand{\V}{\mathcal{V}}
\newcommand{\wOmega}{\tilde\Omega}	

\newcommand{\sharpo}{{\text {\tiny $\sharp$}}}\newcommand{\flatto}{{\text {\tiny$\flat$}}}
\newcommand{\naturallo}{{\text {\tiny$\natural$}}}


\renewcommand{\separe}{}



\thispagestyle{empty} 

\phantom{1}

\vspace{0.4cm}
\noindent{\Large 
{\bf {\boldmath{$\Gamma$}}-convergence for high order phase field fracture: }

\vspace{2pt}
\noindent {\bf continuum and isogeometric formulations}

}

\vspace{38pt}

\begin{small}
{\bf M.~Negri}

{Department of Mathematics -  University of Pavia} 

{Via A.~Ferrata 1 - 27100 Pavia - Italy}

{matteo.negri@unipv.it} 

\vspace{36pt}
\noindent {\bf Abstract.} We consider second order phase field functionals in the continuum setting, and their discretization with  isogeometric tensor product B-splines. We prove that these functionals, continuum and discrete, $\Gamma$-converge to a brittle fracture energy, defined in the space $GSBD^2$. In particular, in the isogeometric setting, since the projection operator is not Lagrangian (i.e., interpolatory) a special construction is needed in order to guarantee that recovery sequences take values in $[0,1]$; convergence holds, as expected, if $h = o (\eps)$, being $h$ the size of the physical mesh and $\eps$ the internal length in the phase field energy.


\end{small}

\vspace{0pt}

 \thispagestyle{empty}



\section{Introduction \label{intro}}

For $\eps>0$ and $\eta_\eps > 0$ we consider the phase field functionals \cite{BordenHughesLandisVerhoosel_CMAME14, LiPecoMillanAriasArroyo_IJNME15} 
\begin{equation} \label{e.i1}
	\F_\eps ( u , v) = \int_{\Omega}  ( v^2 + \eta_\eps ) \, W (\strain (u) )  \, dx+  \tfrac14  G_c  \int_{\Omega}  \eps^{-1} |v-1|^2 + 2 \eps | \nabla v|^2 + \eps^3 | \Delta v |^2  \, dx  ,
\end{equation}
where $\Omega \subset \R^3$, $u \in \U = H^1(\Omega , \mathbb{R}^3)$, and $v \in \V = H^2(\Omega ; [0,1] )$. Here $W (\strain) = \tfrac12 \strain : \mathbf{C} \strain$ is a linear elastic energy density (non-necessarily isotropic) while $G_c >0$ is toughness. As a first result, we show that for $\eta_\eps = o (\eps)$ the $\Gamma$-limit \cite{DalMaso93, Braides98} of $\F_\eps$ (as $\eps \to 0$ and with respect to the strong $L^2$-topology) is the brittle fracture energy 
\begin{equation} \label{e.i2}
	F (u) =  \int_{\Omega \setminus J(u)} W (\strain (u) ) \, dx + G_c \H^{2} ( J(u)) , \ \text{ for } u \in GSBD^2 ( \Omega) ,
\end{equation}
where $J(u)$ denotes the jump set of $u$ and $\H^2$ denotes the Hausdorff measure (roughly speaking the area of $J(u)$). Our proof employs a classical approach in $\Gamma$-convergence: the $\Gamma$-liminf inequality is obtained by  slicing \cite{DalMaso_JEMS13}, together with a one dimensional liminf estimate, while the $\Gamma$-limsup inequality is obtained by density \cite{Iurlano_CVPDE12, ChambolleCrismale_ARMA19}, together with a regularization of the one dimensional optimal profile.  We remark that the slicing technique is made possible by the definition itself of $GSBD$ fields \cite{DalMaso_JEMS13} and by a localization argument which allows to employ the full Hessian instead of the Laplacian (upon introducing an arbitrarily small error). 

Our second result is closer to computational fracture propagation, and above all to \cite{BordenHughesLandisVerhoosel_CMAME14, LiPecoMillanAriasArroyo_IJNME15}. We consider the discretizations 
\begin{equation} \label{e.4}
	\F_{\eps, h} ( u_h , v_h  )  = \int_{\Omega}  ( v_h^2 + \eta_\eps ) \, W (\strain (u_h) )  \, dx+ \tfrac14 G_c  \int_{\Omega}  \eps^{-1} |v_h -1|^2 + 2 \eps | \nabla v_h|^2 + \eps^3 | \Delta v_h |^2  \, dx  ,
\end{equation}
obtained by restriction of the functionals $\F_\eps$ to discrete spaces $\U_h \subset H^1(\Omega , \mathbb{R}^3)$ and $\V_h \subset H^2(\Omega ; [0,1] )$ of isogeometric tensor product B-splines, which are very natural and efficient for high order problems. In the discrete setting, we show that for $\eta_\eps=o(\eps)$ and $h=o(\eps)$ (the element size) the $\Gamma$-limit of $\F_{\eps,h}$ is again the above Griffith's functional $F$ in $GSBD^2$. Comparing with the continuum setting, the discrete $\Gamma$-limsup inequality requires to take into account the fact that  ``interpolation" in the space of tensor product B-splines does not preserve $L^\infty$-bounds; as a consequence the projection $v_h$ of the continuum phase-field profile $v$, which is a natural candidate for the recovery sequence, may not take value in $[0,1]$. This technical issue is by-passed using an {\it ad hoc} local modification of $v_h$, at the price of introducing an additional approximation error, vanishing in the limit for $\eps \to 0$. We stress the fact that the condition $h = o(\eps)$ is necessary and natural in order to guarantee a good enough approximation $v_h$ of the field function $v$ in the transition layer, which is indeed of order $\eps$. In the applications this condition is often guaranteed by $h$-adaptive mesh refinement in a neighbourhood of the crack tip, 
while, from a theoretical point of view, it appears also in the finite element approximation \cite{BellettiniCoscia_NFAO94} of the Mumford-Shah functional.
Our result, with minor modifications, holds also for $C^1$ finite elements and, as a by-product, gives and alternative proof of \cite{BellettiniCoscia_NFAO94}. 


\medskip

\separe
%

In a broad view, the $\Gamma$-convergence result in the continuum Sobolev space setting, fits into a prolific line of research, tracing back to \cite{AmbrosioTortorelli_CPAM90} with the approximation in the sense of $\Gamma$-convergence of the Mumford-Shah functional 
\begin{equation} \label{e.i3}
	MS (u) = \tfrac12 \int_{\Omega}  | \nabla u |^2   \, dx +  G_c \mathcal{H}^1 ( J(u) ) , \quad \text{for } u \in GSBV^2 (\Omega) ,
\end{equation}
by means of the (first order) Ambrosio-Tortorelli functional 
\begin{equation}  \label{e.i4} 
 \tfrac12 \int_{\Omega}  ( v^2 + \eta_\eps ) | \nabla u |^2   \, dx +  \tfrac14 G_c \! \int_{\Omega}  \eps^{-1} |v-1|^2 +  \eps | \nabla v|^2 \, dx,  
\end{equation}
where $u \in H^1(\Omega)$  and $v \in H^1(\Omega, [0,1])$. Note that here $u$ is a scalar. From the technical point of view, switching from the scalar Mumford-Shah functional \eqref{e.i3} to its vectorial counterpart \eqref{e.i2} is not as simple as it may seem. Indeed, a complete $\Gamma$-convergence proof for (first order) vectorial phase field energies of the form 
\begin{equation} \label{e.i5}
	\int_{\Omega}  ( v^2 + \eta_\eps ) \, W (\strain (u) )  \, dx +  \tfrac14  \, G_c \! \int_{\Omega}  \eps^{-1} |v-1|^2 + \eps | \nabla v|^2 \,dx  ,
\end{equation}
was obtained several years after \cite{AmbrosioTortorelli_CPAM90}, first by \cite{Chambolle_JMPA04} in the framework of the space $SBD^2$ and later by \cite{DalMasoIurlano_CPAA13, ChambolleCrismale_ARMA19} in the framework of the larger space $GSBD^2$ (after $GSBD$ was introduced in \cite{DalMaso_JEMS13}). 
Along this line of research, a further important result for fracture has been obtained in \cite{CCF} considering the energy
$$
	\tfrac12 \int_\Omega (v^2 + \eta_\eps) \big( \mu |\strain_{d} (u)|^{2}+\kappa|\strain_{v}^{+}(u)|^{2} \big)  + \kappa|\strain_{v}^{-} (u)|^{2} \, dx  +  \tfrac14  \, G_c \! \int_{\Omega}  \eps^{-1} |v-1|^2 + \eps | \nabla v|^2 \,dx  ,
$$
where $\strain_{v}(u):= \tfrac{1}{2} \mathrm{tr} \strain(u) I$ and $\strain_{d}(u):= \strain(u) - \strain_{v}(u)$ give the volumetric and deviatoric splitting of the strain, while $\strain_{v}^{\pm}(u)$ denotes the positive and negative part. In this case the $\Gamma$-limit  \cite{CCF} takes the form 
\begin{displaymath}
 \tfrac{1}{2}\int_{\Om} \big( 2\mu|\strain_{d}(u)|^{2} + \kappa |\strain_{v}(u)|^{2}\big) \, dx + G_c  \mathcal{H}^{1}(J (u)), \quad  \text{if $ (u^+ - u^- ) \cdot \nu_{u} \geq 0$ in~$J(u)$.}
\end{displaymath}
The constraint $( u^+ - u^- ) \cdot \nu_{u} \ge 0$ provides an (infinitesimal)  non-interpenetration condition on the crack faces, in terms of the crack opening $ (u^+ - u^- )$ along the normal $\nu_{u}$ to $J(u)$.
For difficult technical reasons this $\Gamma$-convergence result holds for $\Omega \subset \R^2$ and for displacement fields in $SBD^2 \cap L^\infty$. 

\separe

As far as second order phase-field functionals, the literature is not as rich as that concerned with first order. One of first results is that of   \cite{FonsecaMantegazza_SJMA00}, dealing with the convergence of second order Modica-Mortola energies, in $H^2(\Omega, [0,1])$, to the perimeter functional, in $BV(\Omega,\{0,1\})$. For the Mumford-Shah functional \eqref{e.i3}  a $\Gamma$-convergence proof for the second order functionals 
\begin{gather*}
\tfrac12 \int_{\Omega}  v^2 | \nabla u |^2   \, dx +  \tfrac14 G_c \! \int_{\Omega}  \eps^{-1} |v-1|^2 +  \eps^{3} | \Delta v|^2 \, dx , \\
\tfrac12 \int_{\Omega}  v^2 | \nabla u |^2   \, dx +  \tfrac14 G_c \! \int_{\Omega}  \eps^{-1} |v-1|^2 +  \eps^{3} | D^2 v|^2 \, dx ,
\end{gather*}
for $u \in H^1(\Omega)$ and  $v \in H^2(\Omega)$ has been proven in \cite{BurgerEspositoZeppieri_MMS15}; note that in these cases the phase-field function $v$ is not constrained to take values in $[0,1]$.  A more general result, for a wider class of free-discontinuity problems, has been recently obtained in \cite{Bach_ESAIMCOCV18}.


\separe

%


\medskip


%
The interest for phase-field functionals is strictly related to the applications. Initially, energies like \eqref{e.i4} have been used in image segmentation problems, e.g.~\cite{March_VC92},
later, after \cite{BourdFrancMar00}, they spread in fracture mechanics, see e.g.,~\cite{MieheWelschingerHofacker10, KuhnMueller_EFM10, SicsicMarigoMaurini_JMPS14, BordenHughesLandisVerhoosel_CMAME14, LiPecoMillanAriasArroyo_IJNME15, MikelicWheelerWick_N15, KiendlAmbatiDeLorenzisGomezReali_CMAME16},  the book \cite{BourdFrancMar08} and the review \cite{AmbatiGerasimovDeLorenzi_CM15}. In this perspective, $\Gamma$-convergence provides a rigorous mathematical framework to prove that phase-field energies are consistent with sharp-crack (free discontinuity) energies. On the other hand, applications in fracture mechanics require, beside energy, an evolution which governs the propagation of the crack. For phase field fracture, this is usually obtained by (time discrete) incremental  problems, based on alternate minimization, or staggered, schemes \cite{BourdFrancMar00}. A characterization of the time-continuous evolution (in the limit as the time step vanishes) has been proved for first order phase-field functionals in \cite{LarsenOrtnerSuli_M3AS10} (for the dynamic case) and in \cite{N_LNACM16, KneesNegri_15, AlmiNegri_19} (for the quasi-static case). Finally, we remark that algorithms based on second order functionals proved to be numerically very efficient; indeed, alternate minimization schemes converge to an equilibrium configuration faster than first order problems (see e.g., \cite[Figure~10 and Tables~4, 5]{BordenHughesLandisVerhoosel_CMAME14} and similarly \cite[Table~1]{BurgerEspositoZeppieri_MMS15}). 





%



\tableofcontents


\section{Setting and statement of the {\boldmath{$\Gamma$}}-convergence results \label{setting}}

\subsection{Continuum setting} 

We assume that the reference domain $\Omega \subset \R^3$ is open, bounded and connected. 
%
The space of admissible {\it continuum displacements} is given by $\U = H^1 (\Omega, \R^3) $ while the space of admissible {\it phase-field functions} is $ \V =  H^2 ( \Omega) $. Note that functions in $\V$ do not necessarily satisfy the constraint $0 \le v \le 1$, which is taken into account in the functional \eqref{e.Fepstilde}. 
The space of admissible {\it discontinuous displacements} is instead provided by  $GSBD^2 ( \Omega) $ 
(see Appendix  \ref{A}, for the definition and the basic properties of this space, and \cite{DalMaso_JEMS13} for the original work).

\separe

For technical reasons, natural in $\Gamma$-convergence, we will employ the ``extended" functionals $\F_\eps $ and $\F$ defined in $L^2(\Omega , \R^3) \times L^2(\Omega)$ and given by 
\begin{align} \label{e.Fepstilde}
	\F_\eps ( u , v) = 
	{\displaystyle \int_{\Omega}  ( v^2 + \eta_\eps  ) \, W (\strain (u) )  \, dx +  \int_{\Omega}  \eps^{-1} |v-1|^2 + 2 \eps | \nabla v|^2 + \eps^3 | \Delta v |^2  \, dx } 
\end{align}
if $(u,v) \in \U \times \V$ and $0 \le v \le 1$, and $\F_\eps ( u , v) = + \infty$ otherwise; 
\begin{align} \label{e.Flim}
	\F ( u , v)  = \begin{cases} 
	 {\displaystyle \int_{\Omega \setminus J(u)} \hspace{-6pt} W (\strain (u) ) \, dx + 4 \hspace{0.3pt} \H^2 ( J(u)) } & \text{if $u \in GSBD^2 ( \Omega)$ and $v =1$ a.e.~in $\Omega$}, \\
	+\infty	& \text{otherwise.}
	\end{cases}
	\
\end{align}
We will assume that $W$ is coercive and continuous in $\R^{3 \times 3}_{sym}$, i.e.~that $c_1 | E |^2 \le W  (E) \le c_2 | E |^2$ for $c_i >0$ and for every $E \in \R^{3 \times 3}_{sym}$.

\begin{remark} \normalfont The choice of $L^2$ in the definition of $\F_\eps$ and $\F$ is due to the fact that $\Gamma$-convergence will be proven with respect to the $L^2$-norm, which seems general enough for our applications. More general choices are also feasible: for instance, taking full advantage of the generality of $GSBD$ spaces, the functionals $\F_\eps$ and $\F$ could be defined in the metric space of measurable vector fields endowed with convergence in measure \cite{ChambolleCrismale_ARMA19}.
\end{remark}


Our main result in the continuum setting is stated in the following Theorem.

\begin{theorem} \label{t.THEO} For $\eta_\eps = o (\eps)$ the functionals $\F_\eps$ $\Gamma$-converge to $\F$ (as $\eps \to 0$) with respect to the strong topology of $L^2 ( \Omega , \R^3) \times L^2( \Omega)$.
\end{theorem}

\begin{remark} Analogous convergence results hold for $\Omega \subset \R^N$ for $N=1,2$, with volume loads in $L^2$ and with Dirichlet boundary conditions for the displacement field \cite{ChambolleCrismale_ARMA19}.
As a standard by-product of $\Gamma$-convergence we have, upon compactness, the strong convergence of minimizer.
\end{remark}

\subsection{Isogeometric quadratic tensor product B-splines}



We follow the assumptions and notation of \cite{BazilevsBeiraoCottrellHughesSangalli_M3AS06} (see also \cite{BeiraoBuffaSangalliVazquez_16}). Let $(0,1)^3$ be the (parametrizing) patch  and let $\mathcal{Q}_h =\{ Q \}$ be a family of uniformly shape regular  meshes of elements $Q$ with diameter $h_Q \le h$; shape regularity means that  the ratio between the length of the edges and the diameter is bounded (from below) uniformly with respect to $Q$ and $\mathcal{Q}_h$.
Let $\mathbf{F}: (0,1)^3 \to \Omega$ be the parametrization map for the physical domain $\Omega$ and denote by $K = \mathbf{F} (Q)$ the elements of the physical mesh $\mathcal K_h = \{ K \}$.  We assume that globally  (from $(0,1)^3$ to $\Omega$) the map $\mathbf{F}$ is 
a diffeomorphism of class $W^{2,\infty}$. As a consequence the family $\mathcal K_h$ is still shape regular and $h_K \le C h$ uniformly with respect to $K$ and $\mathcal K_h$.

We will not enter into the details about the generation of the spaces of quadratic (tensor product) $B$-splines on $\mathcal K_h$ since it is not crucial for our analysis, the reader will find a brief description in \cite{BazilevsBeiraoCottrellHughesSangalli_M3AS06} and a comprehensive treatise in \cite{Schumaker_07}. We will denote by $\U_h $ and $\V_h$ (on the physical meshes $\mathcal K_h$) the discrete spaces of $B$-splines for the displacement field and the phase-field function respectively.  It is important to remark that, in general, functions $v_h \in \V_h$ are allowed to take any real value and thus they may not satisfy the constraint $0 \le v_h \le 1$, which will be imposed in the functional \eqref{e.FhFE}.

We denote by $\tilde{K} \subset \Omega$ the extended support of $K \in \mathcal{K}_h$, i.e.~the union of the supports of the basis functions (of both $\U_h$ and $\V_h$) whose support intersects $K$. We remark that $\tilde K \subset \Omega$ and that $\tilde K \subset \{  \mathrm{dist}( x , K) \le \tilde C h  \} $ for $\tilde C >0$ independent of $K$ and $\mathcal{K}_h$. By \cite[Theorem 3.1]{BazilevsBeiraoCottrellHughesSangalli_M3AS06} we know that there exists a linear approximation operator $\Pi_{\hspace{1.2pt}\U_h} : H^2 (\Omega,\R^3) \to \U_h$ such that for every $0 \le k \le l \le 2$  and every element $K$ of $\mathcal K_h$ it holds
\begin{equation} \label{e.err-Uhloc}
	|  u - \Pi_{\hspace{1.2pt}\U_h} u |_{H^k(K,\,\R^3)} \le C h^{l-k} \| u \|_{H^l (\tilde{K}, \,\R^3)} .
\end{equation}
Similarly, there exists a linear approximation operator $\Pi_{\hspace{0pt}\V_h} : H^3 (\Omega) \to \V_h$ such that for every $0 \le k \le l \le 3$  and every element $K$ of $\mathcal K_h$ it holds
\begin{equation} \label{e.err-Vhloc}
	|  v - \Pi_{\hspace{0pt}\V_h} v |_{H^k(K)} \le C h^{l-k} \| v \|_{H^l (\tilde{K})}  .
\end{equation}
Note that in the previous estimates the norms in the right-hand side are evaluated in the extended element $\tilde{K}$.
Clearly, from local estimates we get also the global ones:
\begin{equation} \label{e.UVhOm}
	|  u - \Pi_{\hspace{1.2pt}\U_h} u |_{H^k(\Omega,\,\R^3)} \le C h^{l-k} \| u \|_{H^l (\Omega, \,\R^3)} .
	\qquad
	|  v - \Pi_{\hspace{0pt}\V_h} v |_{H^k(\Omega)} \le C h^{l-k} \| v \|_{H^l (\Omega)}  .
\end{equation}

\begin{remark} \normalfont Note that, even if $v \in H^3(\Omega)$ takes values in $[0,1]$, in general  $\Pi_{\V_h} v$ does not take values in $[0,1]$ even if the basis functions do. Indeed, high order ''interpolation" in spline or polynomial spaces is not Lagrangian (i.e., interpolatory), it is rather a projection operator which in general does not preserve ordering and $L^\infty$-bounds (see for instance \cite[\S 12]{Schumaker_07}).  A similar issue occurs also for $C^1$ finite elements. In \S \ref{fe} we will provide an ``ad hoc" local modification of the projection $\Pi_{\V_h} v$ (for a special function $v$) taking values in $[0,1]$. 
\end{remark}

Since the elements are uniformly ``equivalent'' to a reference element, through the diffeomorphism $\mathbf{F}$, by a simple change of variable and by Sobolev embedding (in a reference element) it is immediate to see that there exists a constant $C>0$ (independent of $h>0$) such that 
\begin{equation} \label{l.SobK}
	\| z \|_{L^\infty(K)} 	\le	C \big( h^{-3} \|  z \|^2_{L^2 (K)} + h^{-1} | z |^2_{H^1 (K)}  
	+ h |  z |^2_{H^2 (K)} \big)^{1/2}
\end{equation}
for  every $K \in \mathcal{K}_h$ and every $z \in H^2(K)$. Note that this estimate holds for every function in $H^2(\Omega)$ and not only for B-splines.

\medskip
At this point we can introduce the discrete functionals $\mathcal{F}_{\eps,h}$ given by 
\begin{equation} \label{e.FhFE}
	\F_{\eps, h} ( u_h , v_h  )  = 
		\int_{\Omega}  ( v_h^2 + \eta ) \, W (\strain (u_h) )  \, dx 
+ \int_{\Omega}  \eps^{-1} |v_h -1|^2 + 2 \eps | \nabla v_h|^2 + \eps^3 | \Delta v_h |^2  \, dx 
\end{equation}
if $(u_h , v_h) \in \U_h \times \V_h$ and $0 \le v_h \le 1$, and by $ \F_{\eps, h} ( u_h , v_h  )  = +\infty$ otherwise in $L^2 (\Omega , \R^3 ) \times L^2 (\Omega)$.
Note that $\F_{\eps,h}$ is just the restriction of the functional $\F_\eps$ to $\U_h \times \V_h$. The convergence result is the following.

\begin{theorem} \label{t.THEOh} If $\eta = o(\eps)$ and $h  = o (\eps)$ the functionals $\F_{\eps, h}$ $\Gamma$-converge  to $\F$ (as $\eps \to 0$) with respect to the strong topology of $L^2 ( \Omega , \R^3) \times L^2( \Omega)$.
\end{theorem}

The proof of the previous Theorem will follow from Proposition \ref{p.Glinf} and Proposition \ref{p.Glsup}.

\begin{remark} \normalfont 
The condition $h= o(\eps)$, which appears also in \cite{BellettiniCoscia_NFAO94}, allows to have an accurate approximation of the transition layer of the phase-field variable; in practice it should be satisfied only in a neighbourhood on the discontinuity set and often is obtained by local  $h$-refinement, e.g. \cite{ArtinaFornasierMichelettiPerotto_SIAMJSC15, BourdinChambolle_NM00, 
BurkOrtnSuel10, PaulZimmermannMandadapuHughesLandisSauer_19}.  \end{remark}

\subsection{Finite Elements}

The proofs contained in \S\,\ref{fe} have been written in the context of isogeometric tensor product B-splines, because this is the setting of \cite{BordenHughesLandisVerhoosel_CMAME14}. 
Actually, a convergence result like Theorem \ref{t.THEOh} holds, as a by-product, also for finite element spaces (roughly speaking, by replacing in the proofs the extended support $\tilde{K}$ with $K$). More precisely, let $\mathcal{K}_h = \{ K \}$ be a regular family of (triangular or quadrilateral) affine equivalent finite elements in the physical domain $\Omega$. Denote again by $\mathcal{U}_h \subset H^1(\Omega , \R^3)$ and by $\mathcal{V}_h \subset H^2(\Omega)$ the finite element spaces for the displacement fields and phase field functions respectively.  We assume also that there exists a linear approximation operator $\Pi_{\hspace{1.2pt}\U_h} : H^2 (\Omega,\R^3) \to \U_h$ such that for every $0 \le k \le l \le 2$  and every element $K$ of $\mathcal K_h$ it holds
\begin{equation} \label{e.err-Uhloc-FE}
	|  u - \Pi_{\hspace{1.2pt}\U_h} u |_{H^k(K,\,\R^3)} \le C h^{l-k} \| u \|_{H^l (K, \,\R^3)} 
\end{equation}
and that there exists a linear approximation operator $\Pi_{\hspace{0pt}\V_h} : H^3 (\Omega) \to \V_h$ such that for every $0 \le k \le l \le 3$  and every element $K$ of $\mathcal K_h$ it holds
\begin{equation} \label{e.err-Vhloc-FE}
	|  v - \Pi_{\hspace{0pt}\V_h} v |_{H^k(K)} \le C h^{l-k} \| v \|_{H^l (K)}  .
\end{equation}
We remark that the condition $\V_h \subset H^2(\Omega)$ requires continuity of the gradient across element boundaries, i.e.~$C^1$ finite elements; 
we refer to the classic book \cite{Ciarl78} for some examples of elements, for forth order elliptic problems, enjoying this property together with the previous interpolation estimates. Once again, these elements are not Lagrangian and thus interpolation does not preserve, in general, $L^\infty$-bounds.

As in \eqref{e.FhFE} the discrete functionals $\mathcal{F}_{\eps,h}$ are defined by 
\begin{equation*} 
	\F_{\eps, h} ( u_h , v_h  )  = \int_{\Omega}  ( v_h^2 + \eta ) \, W (\strain (u_h) )  \, dx+  \int_{\Omega}  \eps^{-1} |v_h -1|^2 + 2 \eps | \nabla v_h|^2 + \eps^3 | \Delta v_h |^2  \, dx  
\end{equation*}
if $(u_h , v_h) \in \U_h \times \V_h$ and $0 \le v_h \le 1$, and by $ \F_{\eps, h} ( u_h , v_h  )  = +\infty$ otherwise in $L^2 (\Omega , \R^3 ) \times L^2 (\Omega)$.
Note that $\F_{\eps,h}$ is again the restriction of the functional $\F_\eps$ to $\U_h \times \V_h$.

\begin{theorem} \label{t.THEOhFE} If $\eta = o(\eps)$ and $h  = o (\eps)$ the functionals $\F_{\eps, h}$ $\Gamma$-converge  to $\F$ (as $\eps \to 0$) with respect to the strong topology of $L^2 ( \Omega , \R^3) \times L^2( \Omega)$.
\end{theorem}


\section{Preliminary one dimensional estimates \label{1D}}

For $R \in (0, +\infty]$ consider the functionals $\J_R : H^2 (0,R) \to [0, + \infty)$ given by
\begin{equation}\label{e.J}
	\J_R (w) = \int_{(0,R)}  w^2 + 2 | w' |^2 + | w''|^2 \, dr .
\end{equation}

\begin{lemma} \label{l.opti} Let $w_\infty (r) = e^{-r} ( 1 + r)$. Then 
\begin{equation}\label{e.opti}
   w_\infty \in \argmin \{  \J_\infty (w) :  w (0) = 1, \, w'(0) = 0 \}   \quad \text{and} \quad  \J_\infty ( w_\infty) =2  .
\end{equation} 
\end{lemma}

\proof  The Euler-Lagrange equation for $\J_\infty$ reads $w^{(4)} - 2 w'' + w =0$ whose solutions are of the form 
$ w(r) = e^{r} ( C_1 + C_2 r) + e^{-r} ( C_3 + C_4 r)$. Considering the boundary conditions, the unique solution in $H^2 (0,+\infty)$  is given 
by $w_\infty$. An explicit computation gives $\J_\infty ( w_\infty) =2$. \qed 

Note that $w_\infty$ belongs to $W^{m,\infty} (0,+\infty) \cap C^\infty (0,+\infty)$, for $m$ arbitrarily large, and that $w_\infty$ is monotone decreasing with $\lim_{r \to +\infty} w_\infty (r) = 0$, in particular $0 \le w_\infty \le 1$. 

\bigskip
\begin{figure}[h!]
\includegraphics{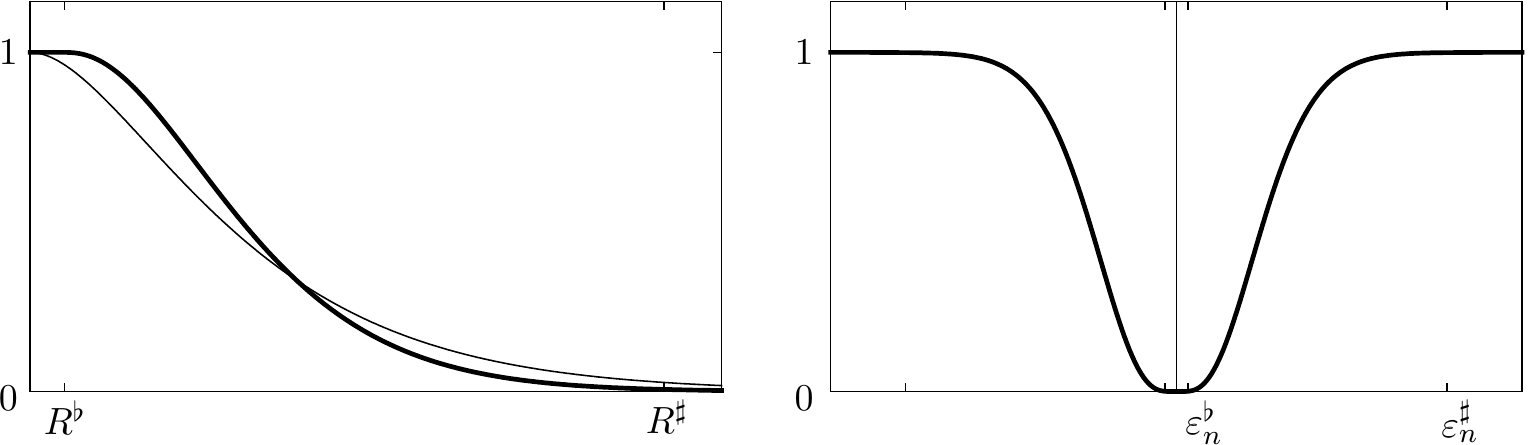}
\caption{\label{F.1} Left: profile of the functions $w_\infty$  from Lemma \ref{l.opti}  and $w$ (solid) from Lemma \ref{l.approx}. Right: profile of a function $z_n$ from Lemma \ref{l.glsup}.}
\end{figure}

The next two lemmas will be used respectively for the $\Gamma$-limsup estimate (Lemma \ref{l.glsup}) and for the $\Gamma$-liminf estimate (Lemma \ref{l.glinf}).

\begin{lemma} \label{l.approx} For $\delta >0 $ there exists $w \in W^{m,\infty} (0,+\infty) \cap C^\infty (0,+\infty)$, for $m$ arbitrarily large, with $0 \le w \le 1$, $w=1$ in $(0, R^\flat)$ and $w=0$ in $(R^\sharp , +\infty)$, for $0 < R^\flat < 1 < R^\sharp$, and such that $2 \le \J_{R^\sharp} (w) = \J_\infty ( w ) < 2 +\delta $. 
\end{lemma}

\proof 
Let $\phi$ be a smooth function in the real line with $\phi (r) = 1 $ for $r < -1 $, $\phi (r) =0 $ for $r>0$ and $0 \le \phi \le 1$. For $0< r_k \to 0^+$ and $R_k \to  +\infty$ let 
\begin{equation} \label{e.wk}
	w_k (r) = \begin{cases}	
	1 & r \le r_k \\
	w_\infty ( r - r_k) \, \phi ( r - R_k) & \text{otherwise.}
	\end{cases}
\end{equation}
Note that $w_k \in H^2 (0,+\infty) \cap C^1 (0,+\infty)$ because  $w_\infty (r) \sim 1 - r^2$ for $r \sim 0$. 
Since $w_k=0$ in $(R_k , +\infty)$ it is an admissible competitor in \eqref{e.opti}, thus we have $\J_{R_k} ( w_k) =  \J_\infty ( w_k ) \ge \J_\infty ( w_\infty) = 2$. It is easy to check that $w_k \to w_\infty$ strongly in $H^2(0,+\infty)$ and thus $\J_\infty ( w_k ) \to \J_\infty ( w_\infty)$, by continuity of $\J_\infty$.

For $0 < s_n \to 0^+$ let $\rho_n (r) = \rho ( r / s_n) / s_n$ be a sequence of smooth mollifiers. Denote $w_{k,n} = w_k * \rho_n$. Clearly $w_{k,n} \to w_k$ in $H^2(0,+\infty)$ and thus $\J_\infty ( w_{k,n} ) \to \J_\infty (w_k)$. Moreover, 
$w_{k,n}' = w'_k * \rho_n$ is continuous with compact support. As a consequence  $w_{k,n} \in W^{1,\infty} (0,+\infty)$. The same argument holds for the derivative of any order, hence $w_{n,k} \in W^{m,\infty} (0,+\infty) \cap C^\infty (0,+\infty)$, for $m$ arbitrarily large.

It is then sufficient to choose $w= w_{k,n}$ for $ k $ and $ n $ sufficiently large.  \qed 

\begin{lemma} \label{l.tilde} 

For $R_n \to +\infty$  let $w_n \in H^2(0,R_n)$ such that 
\begin{align*}
	\lim_{n \to +\infty} w_n (0) =1 , \quad  w'_n (0) = 0 , \quad  \lim_{n \to +\infty} w_n (R_n) = 0 , \quad  \lim_{n \to +\infty} w'_n (R_n) = 0 . 
\end{align*}
Then $\liminf_{n \to +\infty} \J_{R_n} ( w_n) \ge \J_\infty ( w_\infty) = 2$. 
\end{lemma}

\proof  By classical results on Sobolev functions, there exists $C>0$ and a lifting $z_n \in H^2 (0,+\infty)$ with 
$z_n(0) = w_n(R_n)$, $ z'_n(0) = w'_n(R_n)$ and $\| z_n \|_{H^2 (0,+\infty)} \le C ( \hspace{0.5pt}  | w_n(R_n) | + | w'_n(R_n)| \hspace{0.5pt} )  . $
Hence
$$
	Z_n  = \int_{(0,+\infty)}  z_n^2 + 2 | z_n' |^2 + | z_n''|^2 \, dr  \to 0 .
$$
Let $\tilde{w}_n \in H^2(0,+\infty)$ be the extension of $w_n$ given by $\tilde w_n ( r ) = z_n ( r - R_n) $ for $r \in (R_n, +\infty)$.
Denote $\lambda_n =  1 / w_n(0) $. Clearly $\lambda_n \tilde w_n(0) =1$ and $\lambda_n \tilde w'_n (0) =0$, moreover  
$$
	\J_\infty ( \lambda_n \tilde w_n) = \lambda_n^2 \J_\infty ( \tilde w_n) = \lambda_n^2 ( \J_{R_n} ( w_n) + Z_n ) .
$$
As $\lambda_n \to 1$ and $Z_n \to 0$ by minimality of $w_\infty$ we have
$$\liminf_{n \to +\infty} \J_{R_n} (w_n) = \liminf_{n \to +\infty}\lambda_n^{-2} \J_{\infty} ( \tilde w_n) - Z_n  \ge \liminf_{n \to +\infty} \J_\infty ( \lambda_n \tilde w_n) \ge  \J_\infty ( w_\infty) , $$ 
which concludes the proof. \qed

\subsection{An approximate limsup inequality}

\begin{lemma} \label{l.glsup}  For $\delta > 0$ let $w$ be the function provided by Lemma \ref{l.approx}. 
Define $z_n (s) = 1 - w (|s|/\eps_n)$. Then $z_n \in C^\infty ( \R, [0,1])$, $z_n \to 1$ in $L^2_{loc} (\mathbb{R})$ and 
$$
   \lim_{n \to +\infty}  \int_{\R}  \eps_n^{-1} |z_n-1|^2 + 2 \eps_n | z_n'|^2 + \eps_n^3 | z_n'' |^2 \, ds  < 4 + 2 \delta .
$$
Moreover, there exists $C>0$ such that $ \| z_n^{(k)} \|_\infty \le C \eps_n^k$ for $0 \le k \le 3$. 
Finally, denoting $\eps_n^\flat = \eps_n R^\flat$ and $\eps_n^\sharp = \eps_n R^\sharp$, we have $z_n (s) = 0$ for $ |s| \le \eps_n^\flat$ and $z_n (s) = 1$ for $|s| \ge \eps_n^\sharp$.
\end{lemma} 

\proof Since $w=0$ in $(R^\sharp, +\infty)$ it follows that $z_n (s) =1$ for $|s| > \eps_n R^\sharp = \eps_n^\sharp$. Hence $z_n \to 1$ in $L^2_{loc} (\mathbb{R})$. Since $w=1$ in $(0,R^\flat)$ and $w \in C^\infty (0,+\infty)$ we have $z_n (s) = 0$ for $ |s| \le  \eps_n R^\flat = \eps_n^\flat$ and  $z_n \in C^\infty(\R)$.
The change of variable $s = \eps_n r$ yields
$$
	\int_{(0, +\infty)} \eps_n^{-1} |z_n-1|^2 + 2 \eps_n | z_n'|^2 + \eps_n^3 | z_n'' |^2 \, ds 
	=
	\int_{(0,+\infty )}  | w |^2 + 2 | w'|^2 +  | w'' |^2 \, dr  = \J_\infty (w) < 2 + \delta
$$
%
which provides the first estimate.  The estimate for the derivatives can be derived in a similar way by a change of variable.  \qed

%

\subsection{
A liminf inequality}


Let $I=(a,b)$, with $a,b \in \mathbb{R}$, and let $\I_\eps : L^2 (I) \times L^2 (I) \to [0,+\infty]$ be defined  by 
\begin{equation*} 
	\I_\eps ( u , z ) = 
                 \begin{cases} 
                {\displaystyle  \tfrac12 \int_{I}  ( z^2 + \eta_\eps) | u' |^2 \, ds +  \int_{I}  \eps^{-1} |z-1|^2 + 2 \eps | z'|^2 + \eps^3 | z'' |^2 \, ds  }
                 &  \text{if $ (u, z) \in H^{1} (I) \times H^2 ( I)$} \\
		+\infty & \text{otherwise.} 
                 \end{cases}
\end{equation*}
Considering a sequence $\eps_n \to 0$ we will denote $\I_n = \I_{\eps_n}$.
%
%

\begin{lemma} \label{l.glinf} If $(u_n , z_n ) \to (u,z)$ in $L^2 (I) \times L^2 (I)$ and $\liminf_{n \to +\infty} \I_n (u_n , z_n ) < +\infty$ then $ u \in SBV^2 (I)$, $0 \le z \le 1$ and 
\begin{equation}\label{e.glinf}
	4 \, \# ( J(u))  \le \liminf_{n \to +\infty} \int_{I}  \eps_n^{-1} |z_n -1|^2 + 2 \eps_n | z_n'|^2 + \eps_n^3 | z_n'' |^2 \, ds  .  
\end{equation}
\end{lemma}

\proof Neglecting the term $| z''_n|^2$ we get 
$$
	\I_n  (u_n , z_n )  \ge \tfrac12 \int_{I}  ( z_n^2 + \eta_n) | u'_n |^2 \, ds + \int_{I} \eps_n^{-1} |z_n-1|^2 + \eps_n | z_n'|^2 \, ds = AT_n  (u_n , z_n )  ,
$$
where the right hand side is a one dimensional Ambrosio-Tortorelli \cite{AmbrosioTortorelli_CPAM90} functional. Invoking for instance \cite[Theorem 3.15]{Braides98} we get that $z=1$ a.e.~in $I$ and that $u \in SBV^2(I)$ with $ \# ( J(u)) <+\infty$.  
Let $J (u) = \{ s_j \}$. For $\delta \ll 1$ consider the disjoint intervals $I_j^\delta = (s_j-\delta , s_j +\delta) \subset (a,b)$. 
Writing
$$
	\I_n (u_n , z_n) \ge 
         \sum_j \int_{I^\delta_j}  \eps_n^{-1} |z_n-1|^2 + 2 \eps_n | z_n'|^2 + \eps_n^3 | z_n'' |^2 \, ds 
$$
we will check that 
\begin{equation}
	 \liminf_{n \to +\infty} \int_{I^\delta_j}  \eps_n^{-1} |z_n-1|^2 + 2 \eps_n | z_n'|^2 + \eps_n^3 | z_n'' |^2 \, ds  \ge 4 \label{e.glinf2}
	\quad \text{for every $j$},
\end{equation}
from which \eqref{e.glinf} follows. 
As a preliminary step, we extract a subsequence (not relabelled) such that $z_n \to 1$ a.e.~in $I$ and such that each  liminf in \eqref{e.glinf2} is actually a limit. 

Fix an interval $I_j^\delta =  ( s_j-\delta , s_j +\delta)$. Assume, without loss of generality, that $s_j =0$ and denote $I_j^\delta = [-\delta , \delta]$. Fix $\delta^  {\flat}$ and $\delta^  {\sharp}$ (independent of $n$) with $0 < \delta^  {\flat} < \delta^  {\sharp} < \delta$ such that $z_n ( \pm \delta^\flat ) \to 1$ and $z_n ( \pm \delta^\sharp ) \to 1$. First, we show that there exist a subsequence (not relabelled) and for every $n \in \mathbb{N}$ a couple of points, $s^  {\flat}_n \in ( -\delta^  {\flat} , \delta^  {\flat})$ and $s^{\sharp}_n \in (\delta^{\sharp}  , \delta )$, such that  
\begin{equation}  \label{e.qws}
       z_n (s^  {\flat}_n) \to 0,  \quad z'_n(s^  {\flat}_n) = 0, \quad 
       z_n (s^  {\sharp}_n) \to 1,   \quad  | z'_n (s^  {\sharp}_n) | \le 1 . 
\end{equation} 

Let $s^  {\flat}_n \in \argmin \{ z_n (s) : s \in [  -\delta^  {\flat} , \delta^  {\flat}]\}$ (which exists by continuity of $z_n$). Let us see that $z_n (s^  {\flat}_n) \to 0$.
Assume by contradiction that there exists a subsequence (not relabelled) such that $\min \{ z_n (s) : s \in  [ -\delta^  {\flat} , \delta^  {\flat}] \} \ge  C >0$ for every $n \in \mathbb{N}$. Then 
$$
		\I_n ( u_n , z_n) \ge \tfrac12 \int_{ ( -\delta^{\flat} , \delta^  {\flat})}    (z^2_n + \eta_n) | u'_n|^2 \, ds  \ge C  \int_{  (-\delta^  {\flat} , \delta^  {\flat})}  | u'_n|^2 \, ds .
$$
Since $ \I_n ( u_n , z_n) $ is bounded it follows that $u_n$ is bounded in $H^1( -\delta^  {\flat} , \delta^  {\flat})$. As consequence its limit $u$ belongs to $H^1( -\delta^  {\flat} , \delta^  {\flat})$, which contradicts the fact that $0 = s_j \in J(u)$.  
Since $z_n ( \pm \delta^  {\flat} ) \to 1$ the minimizer $s^  {\flat}_n$ belongs to the open interval $(-\delta^  {\flat} , \delta^  {\flat})$, thus $z'_n (s^  {\flat}_n) =0$. 

\separe


Let us find $s^\sharp_n$. For $\tau >0$ consider the open set
$E_n^\tau = \{  s \in [\delta^\sharp , \delta] : 1 - \tau <  z_n (s) \}$. It is not restrictive to consider that $E^\tau_n \neq \emptyset$, indeed, since $z_n \to 1$ in measure we have $| [ \delta^\sharp, \delta] \setminus E^\tau_n | \to 0$.
We want to show that there exists a point $s^\sharp_n \in E^\tau_n$ with $| z'_n (s^  {\sharp}_n) | \le 1$.
Assume by contradiction that $z'_n >1$ in $E_n^\tau$; note that under this assumption $z_n$ is monotone and thus  $E^\tau_n$ is connected. Then, choose $s^* \in E^\tau_n$ with $z_n (s^*) \to 1$; since we have $z_n(s) \ge z_n (s^*) + (s - s^*)$ for $s > s^*$, the upper bound $z_n \le 1$ would be violated if $(s- s^*)>0$ because $z_n (s^*) \to 1$. An analogous argument applies ``backwards'' for $z'_n < -1$. In all the other cases, by the continuity of $z'_n$, there exists a point $s_n^\sharp$ in $E^\tau_n$ with $| z'_n (s_n^\sharp) | \le 1$.
%
%
Choosing $\tau_n \to 0^+$ provides the required sequence.

\separe

Define the rescaled functions $w_n (r) =  1 - z_n ( \eps_n  r + s^  {\flat}_n) $ and let $R_n = ( s^  {\sharp}_n - s^  {\flat}_n) / \eps_n  \ge (\delta^  {\sharp}- \delta^  {\flat}) / \eps_n$. Then $R_n \to +\infty$ and by \eqref{e.qws}
\begin{align*}
	w_n(0) = ( 1 - z_n ( s^  {\flat}_n) ) \to 1 , & \quad w'_n(0) = - \eps_n z'_n(s^  {\flat}_n) = 0, \\
\quad w_n(R_n) = ( 1 - z_n (s^  {\sharp}_n) ) \to 0 , & \quad w'_n(R_n) = - \eps_n z'_n(s^  {\sharp}_n) \to 0 .
\end{align*}
By the change of variable $s= \eps_n r + s_n^\flat$ we have 
$$
	\int_{(s^{\flat}_n , \, s^{\sharp}_n)} \eps_n^{-1} |z_n-1|^2 + 2 \eps_n | z_n'|^2 + \eps_n^3 | z_n'' |^2 \, ds 
	=
	\int_{(0,\, R_n)}  |w_n |^2 + 2 | w_n'|^2 +  | w_n'' |^2 \, dr = \J_{R_n} (w_n) .
$$
Invoking Lemma \ref{l.tilde} we get  $\liminf_{n \to +\infty} \J_{R_n} (w_n) \ge 2$. 
%
%
By symmetry,  \eqref{e.glinf2} is proved. \qed



\section{\boldmath{$\Gamma$}-liminf inequality \label{Ginf}}

Consider a sequence $\eps_n \to 0^+$ and let $\eta_n = o ( \eps_n)$. For simplicity we will employ the notation $\F_n$ for $\F_{\eps_n}$. The $\Gamma$-liminf inequality is based on slicing and on the following standard property, employed also in \cite{BurgerEspositoZeppieri_MMS15}: if $v \in H^2_0 (\Omega)$ then 
\begin{equation} \label{l.H20} 
	\int_\Omega | \Delta v |^2 \, dx = \int_\Omega | D^2 v |^2 \, dx ,
\end{equation}
where $| \cdot |$ in the right-hand side denotes Frobenius norm. 




\begin{proposition} \label{p.Glinf} Let $(u_n , v_n ) \in \U \times \V$ such that $u_n \to u$ in $L^2 (\Omega , \R^3 )$ and $v_n \to v$ in $L^2(\Omega)$. If $\F_n ( u_n , v_n)$ is uniformly bounded then $v=1$ a.e.~in $\Omega$, $u \in GSBD^2(\Omega)$ and 
\begin{gather}
	\liminf_{n \to +\infty} \, \F_n ( u_n , v_n) \ge \F(u,v) .
\end{gather} 
\end{proposition}

\proof As $\F_n ( u_n , v_n)$ is bounded, we have $0 \le v_n \le 1$. Using the first order bound
$$
        \F_n ( u_n , v_n) \ge \int_\Omega  (v^2_n + \eta_n) W ( \strain(u_n)) \, dx \,  +  \int_\Omega  \eps_n^{-1} | v_n -1 |^2 + \eps_n |\nabla v_n |^2 \, dx 
$$
and then arguing as in \cite[Theorem 4.3]{Iurlano_CVPDE12} we get that $v_n \to 1$ in $L^2(\Omega)$, $u \in GSBD^2 (\Omega)$ and that $(v^2_n + \eta_n)^{1/2} \strain (u_n) \weakto \strain(u)$ in $L^2 ( \Omega , \R^{3 \times 3}_{sym})$; thus
\begin{equation}  \label{e.zzz}
	\liminf_{n \to +\infty}  \int_\Omega (v^2_n + \eta_n)\, W (\strain (u_n)) \, dx  =  
	\liminf_{n \to +\infty}  \int_\Omega \, W \big(  ( v^2_n + \eta_n)^{1/2} \strain (u_n) \big) \, dx  \ge 
	 \int_{\Omega \setminus J(u)}  W (\strain (u) ) \, dx   .
\end{equation}

To get the right bound for the jump we need also the second derivatives of the phase field. To this end, first we replace (locally) the Laplacian with the norm of Hessian, introducing a small error, vanishing in the limit.
Given an open set $A \subset \subset \Omega$ let $\phi \in C^\infty_c (\Omega)$ with $0 \le \phi \le 1$ and $\phi =1$ on $A$. Using Young's inequality $(a+b)^2 \le (1+\delta) \, a^2 + ( 1 + \delta^{-1}) \, b^2$ for $\delta>0$, we can write
$$
	| \Delta (v_n \phi)|^2   = | ( \Delta v_n ) \phi + 2 \nabla v_n \cdot \nabla \phi + v_n ( \Delta \phi ) |^2 \le (1+\delta) | \Delta v_n |^2 \phi^2 + C_\phi (1 + \delta^{-1}) (   | \nabla v_n|^2 + | v_n |^2 ) ,
$$
where $C_\phi$ depends only on $\phi$.  Hence, for $C_{\delta,\phi} = C_\phi (1+\delta)^{-1} / ( 1 + \delta)$, being $ 0 \le v_n \le 1$ we have 
$$
	\int _\Omega  | \Delta v_n |^2 \, dx 
	  \ge \int _\Omega | \Delta v_n|^2 \phi^2 \, dx 
 \ge (1+\delta)^{-1} \int_\Omega | \Delta ( v_n \phi )|^2 \, dx  -  C_{\delta,\phi}  \int_\Omega  | \nabla v_n|^2 + 1  \, dx .
$$
Moreover, by \eqref{l.H20}
$$
	\int_\Omega  | \Delta ( v_n \phi )|^2 \, dx = \int_\Omega | D^2 ( v_n \phi )|^2 \, dx  \ge  \int_A | D^2 ( v_n \phi )|^2 \, dx  =  \int_A  | D^2 v_n |^2 \, dx.
$$
Hence
$$
	\int _\Omega  | \Delta v_n|^2 \, dx  \ge  (1 + \delta)^{-1} \int_A  | D^2 v_n |^2 \, dx - C_{\delta,\phi}  \int_\Omega  | \nabla v_n|^2 + 1  \, dx .
$$
Thus, for every $\delta >0$, we can write 
\begin{align*}
	\F_n ( u_n , v_n)  
	 = & \int_\Omega  (v^2_n + \eta_n) W ( \strain(u_n)) \, dx \, 
	+ \int_\Omega \eps_n^{-1} | v_n -1 |^2 + \eps_n |\nabla v_n |^2 + \eps_n^3 | \Delta v_n |^2 \, dx  \\
	\ge &  \int_\Omega  (v^2_n + \eta_n) W ( \strain(u_n)) \, dx \, +  \\  & + (1+ \delta)^{-1} \int_A \eps_n^{-1} | v_n -1 |^2 + \eps_n |\nabla v_n |^2 +  \eps_n^3 | D^2 v_n |^2 \, dx - \eps_n^3  C_{\delta,\phi}  \int_\Omega  | \nabla v_n|^2 + 1 \, dx  .
\end{align*}
Being 
 $$ 0 \le   \eps_n^2  \int_\Omega \eps_n | \nabla v_n|^2 +  \eps_n \, dx \le \eps_n^2 (  \F_n (u_n, v_n) +\eps_n  | \Omega|)   \to 0 $$ 
for every open set $A \subset \subset \Omega$ and every $\delta >0$ we have 
\begin{align}
	\liminf_{n \to +\infty} \F_n ( u_n , v_n)  \ge & \,  \liminf_{n \to +\infty}  \int_\Omega (v^2_n + \eta_n) W (\strain (u_n)) \, dx   \, +  \nonumber \\ & + (1 + \delta)^{-1} \liminf_{n \to +\infty} \, \int_A  \eps_n^{-1} | v_n -1 |^2 + \eps_n |\nabla v_n |^2 + \eps_n^3 | D^2 v_n |^2 \, dx  \label{e.dmah}. 
\end{align}
Let us check that, for every $\xi \in \mathbb{S}^2$ we have 
$$
	\liminf_{n \to +\infty} \, \int_A  \eps_n^{-1} | v_n -1 |^2 + \eps_n |\nabla v_n |^2 + \eps_n^3 | D^2 v_n |^2 \, dx  \ge
	4 \int_{J^\xi(u) \cap A} \! | \nu \cdot \xi | \, d \H^2 .
$$
We will use the slicing technique, see \S \ref{A}. For $\xi \in \mathbb{S}^2$ and $y \in \xi^\perp$ we denote  $A^\xi_y = \{ s \in \mathbb{R} : y + s \xi  \in A \}$. Accordingly, let $v_y^\xi (s) = v ( y + s \xi)$ and  $ u^\xi_y ( s) = u ( y + s \xi ) \cdot \xi$. Since  $u_n \to u$ in $L^2(\Omega, \R^3)$ and $v_n \to 1$ in $L^2(\Omega)$ then for every $\xi \in \mathbb{S}^2$ and a.e.~$y \in \xi^\perp$ we have  $(u_n)^\xi_y \to u^\xi_y$ and $(v_n)^\xi_y \to 1$ in $L^2(\Omega^\xi_y)$. Note also that $u^\xi_y$ belongs to $SBV (A^\xi_y)$, by Definition \ref{d.GSBD2def}, and that, for a.e.~$y \in \xi^\perp$, we have (a.e.~in $A^\xi_y$)
$$
		| \nabla v_n | \ge | D_\xi v_n | =  |D (v_n)^\xi_y  |\, ,
	\qquad
		| D^2 v_n | \ge | D^2_{\xi\xi} v_n | = | D^2 (v_n)^\xi_y | .
$$
We remark that replacing the Laplacian with the full Hessian allows to get the previous bound on the second derivative of the slice.
Then Fubini's Theorem yields
\begin{align*}
	\int_A  \eps_n^{-1} | v_n -1 |^2  & + \eps_n |\nabla v_n |^2 + \eps_n^3 | D^2 v_n |^2 \, dx  \ge  \\
	& + \int_{\xi^\perp} \Big( \int_{A^\xi_y} \eps_n^{-1} | (v_n)^\xi_y - 1 |^2 + \eps_n | D (v_n)^\xi_y   |^2 + \eps_n^3 | D^2 (v_n)^\xi_y |^2 \, ds 
	\Big) d \H^2 (y) .
\end{align*}
By Lemma \ref{l.glinf} and Theorem \ref{t.GSBDslic} we get (for a.e.~$y \in \xi^\perp$)
$$
	\liminf_{n \to +\infty} \int_{A^\xi_y} \eps_n^{-1} | (v_n)^\xi_y - 1 |^2 + \eps_n | D (v_n)^\xi_y   |^2 + \eps_n^3 | D^2 (v_n)^\xi_y |^2 \, ds 
	\ge 4 \, \# \big( J ( u^\xi_y)  \cap A^\xi_y \big) =  4 \, \# \big( (J^\xi (u) \cap A )^\xi_y  \,   \big) .
$$
Therefore, Fatou's Lemma and Theorem \ref{t.GSBDslic} give
\begin{align}
	\liminf_{n \to +\infty} \int_A  \eps_n^{-1} | v_n -1 |^2   + \eps_n |\nabla v_n |^2 + \eps_n^3 | D^2 v_n |^2 \, dx  & \ge 4 
	\int_{\xi^\perp} \! \# \big( (J^\xi (u) \cap A)^\xi_y \, \big)\,   d \H^2 (y) \nonumber \\ 
	& = 4 \int_{J^\xi(u) \cap A} \! | \nu \cdot \xi | \, d \H^2 \label{e.previ}.
\end{align}
Using (\ref{e.zzz})-(\ref{e.previ}) and taking the supremum with respect to $A \subset \subset \Omega$ we get 
$$
	\liminf_{n \to +\infty} \F_n ( u_n , v_n)  \ge 
	\int_\Omega  W (\strain (u) ) \, dx  + 4 ( 1 + \delta)^{-1} \int_{J^\xi (u)} \! | \nu \cdot \xi | \, d \H^2 , 
	\qquad
	\text{for every $\xi \in \mathbb{S}^2$.}
$$
%
%
To conclude we will employ a supremum of measures argument, see \cite[Proposition 1.16]{Braides98}. Let $B \subset \Omega$ be an open set. Denote 
$$
	\F_n ( u , v \, | B ) = \int_{B}  ( v_n^2 + \eta_n ) \, W (\strain (u_n) )  \, dx+ \int_{B}  \eps_n^{-1} |v_n-1|^2 + 2 \eps_n | \nabla v_n|^2 + \eps_n^3 | \Delta v_n |^2  \, dx  .
$$
Arguing as above, just replacing $\Omega$ with $B$, we get 
$$
	\liminf_{n \to +\infty} \F_n (u_n , v_n \, | B) \ge \int_B  W (\strain (u) ) \, dx  + 4  ( 1 + \delta)^{-1} \int_{J^\xi (u) \, \cap \, B} \! | \nu \cdot \xi | \, d \H^2 ,
	\qquad
	\text{for every $\xi \in \mathbb{S}^2$.}
$$
By Theorem \ref{t.GSBDslic}, for a.e.~$\xi \in \mathbb{S}^2$ we have $\H^2 ( J (u) \setminus J^\xi (u) ) = 0$ and thus
$$
	\liminf_{n \to +\infty} \F_n (u_n , v_n \, | B) \ge \int_B  W (\strain (u) ) \, dx  + 4  ( 1 + \delta)^{-1} \int_{J (u) \, \cap \, B} \! | \nu \cdot \xi | \, d \H^2 ,
$$
for a.e.~$\xi \in \mathbb{S}^2$ and every $B \subset \Omega$.
Note that $\sup_{\xi} | \nu \cdot \xi | = 1$, even if the supremum is taken with respect to a.e.~$\xi \in \mathbb{S}^2$. Therefore, by \cite[Proposition 1.16]{Braides98} and \eqref{e.dmah} we get
$$
	\liminf_{n \to +\infty}  \F_n ( u_n , v_n) \ge   \int_\Omega  W (\strain (u) ) \, dx   + 4 (1+\delta)^{-1} \H^2 (J (u)) ,
$$
which concludes the proof, by arbitrariness of $\delta$. \qed



\section{\boldmath{$\Gamma$}-limsup inequality \label{Gsup}}

By a standard diagonal argument in the theory of $\Gamma$-convergence together with Theorem \ref{t.GSBDdens} it is enough to prove the limsup estimate stated in the next Proposition.

\begin{proposition} \label{p.Glsup} 
Let $J \subset  \Omega$ be a closed $2$-simplex and let $u \in W^{2,\infty} ( \Omega \setminus  J , \,\R^3)$. There exists $C>0$ (depending only on $J$) such that for every $\delta >0$ there exist $u_n \in W^{2,\infty} (\Omega , \R^3)$ and $v_n \in W^{3,\infty} ( \Omega, [0,1])$ such that $u_n \to u$ in $L^2 ( \Omega, \R^3)$, $v_n \to 1$ in $L^2 (\Omega)$ and 
\begin{equation} \label{e.18}
	\limsup_{n \to +\infty} \F_n ( u_n , v_n) \le \int_{\Omega \setminus J} W (\strain (u) ) \, dx + 4 |J | + C \delta ,
\end{equation}
where $| J | = \H^2(J)$ denotes the area of $J$.
\end{proposition} 

\proof {\bf Step 1.} Assume, without loss of generality, that $J \subset \{ x_3 =0 \}$. By abuse of notation, we consider $J \subset \R^2$ and write the simplex as $J \times \{ 0 \}$.  We also assume, without loss of generality, that $0 \in \mathring{J}$ (the interior of $J$ in the topology of $\R^2$). 
For $\delta >0$, denote $J_\delta = (1+\delta)J$ and note that $|J_\delta| \le |J_{2\delta}| 
\le |J| + C _J \delta $ where $C_J$ depends on $J$.

Let $\phi_\delta \in C^\infty_c ( \R^2, [0,1])$ with $\phi_\delta = 1$ in $J_\delta$ and $\phi_\delta =0$  in $\R^2 \setminus J_{2\delta}$.  Let $z_n$ (depending on $\delta$) be the sequence provided by Lemma \ref{l.glsup} and define
$$
	v_n ( x_1, x_2, x_3) = \phi_\delta (x_1, x_2) (z_n (x_3) - 1) + 1 . 
$$ 
We denote
$$
     A_n = J_\delta \times ( - \eps_n^\sharp , \eps_n^\sharp) ,
              \qquad
     B_n = (J_{2 \delta} \setminus J_\delta) \times ( - \eps_n^\sharp , \eps_n^\sharp) .
$$
Let $\delta \ll 1$ and $n \gg 1$, in such a way that $(A_n \cup B_n) \subset \Omega$. Note that $0 \le v_n \le 1$, $v_n (x) = z_n (x_3)$ in $A_n$ (because $\phi_\delta =1$ in $J_\delta$), and $v_n (x) =1$ in $\Omega \setminus (A_n \cup B_n)$ (because $z_n (x_3) =1$ if $| x_3| \ge \eps_n^\sharp$ and $\phi_\delta = 0$ in $\R^2 \setminus J_{2 \delta}$). It follows that $v_n \to 1$ in $L^2(\Omega)$. 

For convenience, let us also introduce the functions $\bar\phi_\delta (x) = \phi_\delta (x_1,x_2) $ and $\bar{z}_n (x) = z_n(x_3)$, so that we can write $v_n = \bar \phi_\delta (\bar{z}_n - 1) +1$. Clearly, we have $\nabla v_n = \nabla \bar \phi_\delta ( \bar z_n -1) + \bar \phi_\delta \nabla \bar z_n$ and 
$$
	\Delta v_n = \Delta \bar \phi_\delta (\bar z_n-1) + 2 \nabla \bar\phi_\delta \cdot \nabla \bar z_n + \Delta \bar z_n \bar \phi_\delta =  \Delta \bar \phi_\delta ( \bar z_n-1) + \Delta \bar z_n \bar \phi_\delta ,
$$
where, in the second identity, we simply used the fact that  
$ \nabla \bar \phi_\delta \cdot \nabla \bar z_n = ( \partial_1  \phi_\delta, \partial_2 \phi_\delta,0) \cdot (0,0, z'_n)   = 0$.

As $v_n (x) = z_n (x_3)$ in $A_n$, by Lemma \ref{l.glsup} we can write 
\begin{align*}
     \limsup_{n \to +\infty} \int_{A_n} \eps_n^{-1} | v_n -1 |^2 & + 2 \eps_n | \nabla v_n |^2 +  \eps_n^3 | \Delta v_n |^2 \, dx  & \\ & \le \limsup_{n \to +\infty } | J _\delta |  \int_\mathbb{R} \eps_n^{-1} | z_n -1 |^2 + 2 \eps_n |  z'_n |^2 +  \eps_n^3 | z''_n |^2 \, dx_3 \le 4 | J | + C \delta ,
\end{align*}
where $C$ depends only on $J$. By definition, 
\begin{gather*}
	| v_n -1 |^2 = | \bar \phi_\delta |^2 | \bar z_n -1 |^2 \le | \bar z_n - 1|^2 , \\
	| \nabla v_n |^2 \le \left( | \nabla \bar \phi_\delta | \, | \bar z_n -1 | +  | \bar \phi_\delta| \, | \nabla \bar z_n |   \right)^2  \le 2 | \nabla \bar \phi_\delta | ^2 + 2 | \nabla \bar z_n | ^2 ,
\end{gather*}
and, in the same way, 
$$
     | \Delta v_n |^2  \le 2 | \Delta \bar \phi_\delta |^2 + 2 | \Delta \bar z_n |^2 . 
$$
Hence, to evaluate the integral over $B_n$ we can estimate separately the terms with $\phi_\delta$ and $z_n$, writing
\begin{align*}
\int_{B_n} \eps_n^{-1} | v_n -1 |^2 & + 2 \eps_n | \nabla v_n |^2 +  \eps_n^3 | \Delta v_n |^2 \, dx  \le \\ 
	& \le 2 | J_{2 \delta} \setminus J_\delta | \int_{\R} \eps_n^{-1} | z_n -1 |^2 + 2 \eps_n | z'_n |^2 +   \eps_n^3 |  z''_n |^2 \, dx_3 \, + \\
& \phantom{\le} \ +  4 \eps_n^\sharp \int_{J_{2\delta} \setminus J_\delta} \eps_n | \nabla \phi_\delta|^2 + \eps_n^3 | \Delta \phi_\delta |^2 \, dx_1 dx_2   .
\end{align*}
As $| J_{2\delta} \setminus J_\delta | \le | J_{2\delta} | - | J_\delta | \le c \delta$, using Lemma \ref{l.glsup} it follows that 
$$
\limsup_{n \to +\infty}  
\int_{B_n} \eps_n^{-1} | v_n -1 |^2+ 2 \eps_n | \nabla v_n |^2 +  \eps_n^3 | \Delta v_n |^2 \, dx  \le C' \delta.
$$
Finally, being $v_n =1$ in $\Omega \setminus (A_n \cup B_n)$ it is obvious that
$$
    \int_{\Omega \setminus (A_n \cup B_n)}  \eps_n^{-1} | v_n -1 |^2  + 2 \eps_n | \nabla v_n |^2 +  \eps_n^3 | \Delta v_n |^2 \, dx = 0 .
$$
In conclusion, 
$$  
    \limsup_{n \to +\infty}  \int_{\Omega}  \eps_n^{-1} | v_n -1 |^2  + 2 \eps_n | \nabla v_n |^2 +  \eps_n^3 | \Delta v_n |^2 \, dx \le 4 | J | + C \delta .
$$

{\bf Step 2.} Let $\psi_\delta \in C^\infty_c ( \R^2 , [0,1])$ with $\psi_\delta = 0$ in $\R^2 \setminus J_{\delta/2}$ and $\psi_\delta =1$ in a neighbourhood of $J$. Let $\xi \in C^\infty ( \R , [0,1])$ with $\xi (s) = 0$ for $| s | \ge 1$ and $\xi(s)=1$ in a neighbourhood of $0$. Denote $\xi_n ( s) = \xi( 2 s / \eps_n^\flat)$ and define 
$$
	\zeta_n (x_1 , x_2 , x_3) = 1 - \psi_\delta (x_1, x_2) \xi_n ( x_3) .
$$
Denote 
$$
      E_n = J_{\delta/2 } \times ( -\eps_n^\flat /2 , \eps_n^\flat/2  ) .
$$
Note that $\zeta_n \in C^\infty ( \Omega, [ 0,1])$, $\zeta_n =0$ in a neighbourhood of $J \times \{ 0\}$ (because $\xi_n=1$ in a neighbourhood of $0$ and  $\psi_\delta=1$ in a neighbourhood of $J$), $\zeta_n =1$ in  $\Omega \setminus E_n$ (because $\xi_n(x_3)=0$ if $| x_3 | \ge \eps_n^\flat / 2$ and $\psi_\delta =0$ in $\Omega \setminus J_{\delta/2}$); in particular $\zeta_n \to 1$ in $L^2(\Omega)$.
We define
$$
     u_ n = \zeta_n u .
$$
Note that $u_n \in W^{2,\infty} ( \Omega, \R^3)$, by the regularity of $u$ and because $\zeta_n =0$ in a neighbourhood of $J \times \{ 0 \}$;  moreover $u_n = u$ in $\Omega \setminus E_n$, then
$$
\int_{ \Omega \setminus E_n} ( v_n^2 +\eta_n) W ( \strain(u_n) ) \, dx  \le (1 + \eta_n) \int_{\Omega \setminus E_n} W ( \strain(u) ) \, dx \  \to \  \int_{\Omega \setminus J \times \{ 0 \} } W ( \strain(u) ) \, dx  .
$$
Note that $v_n = 0$ in $J_\delta \times (-\eps_n^\flat, \eps_n^\flat)$ (because $z_n (x_3) =0$ if $| x_3| \le \eps_n^\flat$ and $\phi_\delta = 1$ in $J_{\delta}$) and thus in $E_n$, hence
$$
\int_{E_n} ( v_n^2 +\eta_n) W ( \strain(u_n) ) \, dx \le  \eta_n C \int_{E_n} | D u_n |^2 \, dx .
$$
For convenience, let $\bar \psi_\delta (x) = \psi_\delta ( x_1 , x_2)$ and  $\bar \xi_n (x) = \xi_n (x_3)$, write $\zeta_n = 1 - \bar\psi_\delta \bar\xi_n$, and then $| \nabla \zeta_n |  \le  | \nabla \bar\psi_\delta | + | \nabla \bar\xi_n |$ (because both $\bar{\psi}_\delta$ and $\bar\xi_n$ take values in $[0,1]$). Moreover, by the regularity if $u$ and $\psi_\delta$ we get 
$$
	 | D u_n | \le  | \nabla \zeta_n | \,  | u |  + | \zeta_n | \, | D u | \le C ( 1 + | \nabla \bar \xi_n |) ,
$$
where $C$ depends on $u$ and $\delta$. Hence, we can write 
\begin{align*}
     \eta_n \int_{E_n} | D u_n |^2 \, dx 
	& \le \eta_n C \int_{E_n} 1 + | \nabla \bar\xi_n |^2 \, dx \le  \eta_n C' | J_{\delta/2}|  \int_0^{\eps_n^\flat/2} 1 + | \xi'_n (x_3) |^2 dx_3  .
\end{align*}
Let us estimate
$$
     \int_0^{\eps_n^\flat/2} | \xi'_n (x_3) |^2 dx_3 = (\eps_n^\flat/2)^{-2} \int_0^{\eps_n^\flat/2} | \xi' (2 x_3 / \eps_n^\flat) |^2  dx_3 = (\eps_n^\flat/2)^{-1} \int_0^1 | \xi' (s) |^2 \, ds \le C (\eps_n^\flat)^{-1} .
$$
As $\eps_n^\flat = R_n^\flat \eps_n$ and $\eta_n / \eps_n \to 0$,  it follows that 
$$
   \limsup_{n \to +\infty} \int_{E_n} ( v_n^2 +\eta_n) W ( \strain(u_n) ) \, dx \le \limsup_{n \to +\infty}  \eta_n C \int_{E_n} | D u_n |^2 \, dx \le C' \eta_n \eps_n^{-1} = 0 ,
$$
which concludes the proof. \qed

\begin{remark} \normalfont The fact that $u \in W^{2,\infty} (  \Omega, \R^3)$ and $v \in W^{3,\infty} ( \Omega, [0,1])$, instead of the more natural $H^{1} (  \Omega, \R^3)$ and $H^2( \Omega, [0,1])$ which would be enough for the $\Gamma$-limsup estimate, will be useful to employ the projection operators in $\U_{h_n}$ and $\V_{h_n}$ in the discrete approximation, see \S\,\ref{fe}, together with the next Corollary.
\end{remark}


\begin{corollary} \label{c.corlsup} Let $v_n$ be as in Proposition \ref{p.Glsup}. 
Then, there exists $C>0$ (depending on $\delta$ but independent of $n$) such that 
%
$\| D^k v_n \|_{L^\infty (\Omega)} \le C  \eps_n^{-k}$ for $0 \le k \le 3$ and $\| v_n -1 \|_{H^3(\Omega)}^2 \le C \eps_n^{-5}$. Moreover, $v_n = 1$ in $\Omega \setminus (A_n \cup B_n)$ where $A_n \cup B_n = J_{2 \delta} \times (- \eps_n^\sharp, \eps_n^\sharp)$ and $v_n = 0$ in $C_n = J_{\delta} \times (- \eps_n^\flat , \eps_n^\flat)$.
\end{corollary}

\proof Remember that $v_n (x) = \phi_\delta(x_1, x_2)  ( z_n (x_3) -1) +1$, by Lemma \ref{l.glsup} we get the bound on the $L^\infty$-norm of the derivatives. 
To estimate the $H^3$-norm it is enough to employ the previous bound, remembering that $v_n-1$ is supported in the set $(A_n \cup B_n) \times (-\eps_n^\sharp, \eps_n^\sharp) = J_{2\delta} \times (-\eps_n^\sharp, \eps_n^\sharp) $, whose measure is of order $\eps_n$.  \qed
 
In a similar way we get the following corollary. 

\begin{corollary} \label{c.corlsupbis} Let $u_n$ be as in Proposition \ref{p.Glsup}. 
Then, there exists $C>0$ (depending on $\delta$ but independent of $n$) such that $\| D^k u_n \|_{L^\infty (\Omega)} \le C \eps_n^{-k}$ for $0 \le k \le 2$. Moreover, $u_n = u$ in $\Omega \setminus E_n$ where $E_n = J_{\delta/2} \times (- \eps_n^\flat/2 , \eps_n^\flat/2)$.
\end{corollary}

\section{\boldmath{$\Gamma$-limit of $\F_{\eps,h}$} \label{fe}}

As $\F_{\eps,h}$ is the restriction of $\F_\eps$ to $\U_h \times \V_h$  the $\Gamma$-liminf inequality for $\F_{\eps,h}$ follows directly from Proposition \ref{p.Glinf}. Moreover, as in the continuum setting, it is enough to prove the following $\Gamma$-limsup inequality.

\begin{proposition} \label{p.Glsupiso} Let $\eps_n \to 0^+$, $\eta_n = o (\eps_n)$ and $h_n = o(\eps_n)$. Let $J \subset \Omega$ be a closed $2$-simplex and let $u \in W^{2,\infty} ( \Omega \setminus\!  J \hspace{0.5pt}, \R^3)$. There exists $C>0$ such that for every $\delta >0$ there exist $u_{h_n} \in \U_{h_n}$ and $v_{h_n} \in \V_{h_n}$, with $ 0 \le v_{h_n} \le 1$, such that $u_{h_n} \to u$ in $L^2 ( \Omega, \R^3)$, $v_{h_n} \to 1$ in $L^2 (\Omega)$ and 
\begin{equation} \label{e.18h}
	\limsup_{n \to +\infty} \F_{\eps_n , h_n} ( u_{h_n} , v_{h_n}) \le(1+ \delta)  \int_{\Omega \setminus J} W (\strain (u) ) \, dx + 4 | J | + C \delta .
\end{equation}
\end{proposition}

\proof  We adopt the assumptions and notation employed in the the proof of Proposition \ref{p.Glsup}. 
Let $u_n$ and $v_n$ be provided by Proposition \ref{p.Glsup}.


{\bf Step 1.} 
Let $\Pi_{\V_{h_n}}$ be the interpolation operator in $H^{3} (\Omega)$ and let $w_{h_n} = \Pi_{\V_{h_n}} v_n$.
Remember that $v_n \in W^{3,\infty} (\Omega)$ and thus the interpolation error estimate \eqref{e.err-Vhloc}  gives
$$ 
	|  v_n - w_{h_n} |_{H^k(K)} \le C h^{3-k} \| v_n \|_{H^3 (\tilde{K})}  \quad \text{for $0 \le k \le 3$.}
$$ 
By Corollary \ref{c.corlsup} for $0 \le l \le 3$  we have $\| D^{l} v_n \|_\infty \le C \eps_n^{-3}$ and then
for $0 \le k \le 2$ 
\begin{equation} \label{e.run}
            |  v_n - w_{h_n} |^2_{H^k (K)} \le C \,  h_n^{6- 2k} \| v_n \|^2_{H^3 (\tilde K)} \le  C' \,  h_n^{6- 2k} | \tilde{K} |  \,\eps_n^{-6} \le  C'' \,  h_n^{3 - 2k} (h_n / \eps_n)^{6} .
\end{equation}

Note that, in general, the inequality $0 \le w_{h_n} \le 1$ may not hold everywhere in $\Omega$; to fix this point let us start with an estimate of the error $\| v_n - w_{h_n} \|_{L^\infty(\Omega)}$. For every element $K \in \mathcal{K}_{h_n}$ in the physical domain,  \eqref{l.SobK} provides
\begin{equation} \label{e.stella}
	\|  v_n -  w_{h_n} \|_{L^\infty (K)} \le C \big( h_n^{-3} \|  v_n -  w_{h_n} \|^2_{L^2 (K)} + h_n^{-1} |  v_n -  w_{h_n} |^2_{H^1 (K)}
	+ h_n |  v_n -  w_{h_n} |^2_{H^2 (K)} \big)^{1/2} .
\end{equation} 
Joining \eqref{e.run} and \eqref{e.stella} it follows that
$\|  v_n -  w_{h_n} \|_{L^\infty (K)} \le C (h_n/\eps_n)^3$.            
As the constant $C$ is independent of the element $K$ the previous estimate becomes
\begin{equation} \label{e.linfty}
        \|  v_n -  w_{h_n} \|_{L^\infty (\Omega)} \le C (h_n/\eps_n)^3 = c_n   .
\end{equation}
Hence $\| v_n -  w_{h_n} \|_{L^\infty (\Omega)} \to 0$ and $- c_n \le w_{h_n} \le 1 + c_n$ in $\Omega$. 
As explained in Remark \ref{r.remo} it is not possible to simply rescale $w_{h_n}$ in a way that it takes values in $[0,1]$; we employ instead the following local construction. As a first step, define

\begin{figure} \begin{center}
\includegraphics{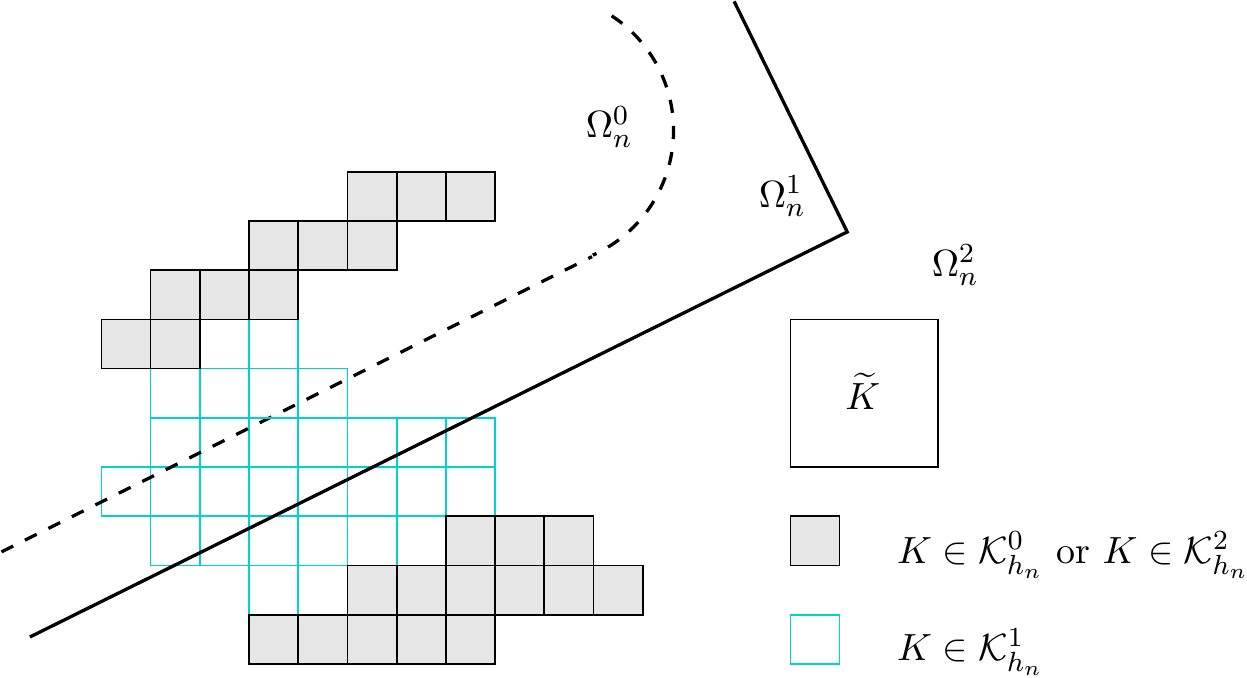}
\end{center} 
\caption{\label{fig3} Sets involved in the proof of Proposition \ref{p.Glsupiso}.} \end{figure}

\begin{gather*}
	\Omega^{0}_n = \{ v_n = 0 \}  , \qquad  \Omega^{1}_n = \{ 0 < v_n < 2 c_n \}  , \qquad  \Omega^{2}_n = \{ 2 c_n \le v_n \le 1 - 2 c_n \},  \\  		\Omega^{3}_n = \{ 1 - 2 c_n  <  v_n < 1  \} , \quad  \Omega^4_n = \{ v_n = 1 \}  .
\end{gather*}
Define also 
\begin{gather*}
	\mathcal{K}^i_{h_n} = \{ K : \tilde{K} \subset \Omega^i_n  \}  \quad \text{ for $i=0,2,4$,}  \qquad 
	\mathcal{K}^i_{h_n} = \{ K :  \tilde K \cap \Omega^i_n \ne \emptyset \}  \quad \text{ for $i=1,3$}. 
\end{gather*}
Note that the previous definitions depends on the extended elements $\tilde K$. First, we check that the families $\mathcal{K}^i_{h_n}$provide a disjoint partition of $\mathcal{K}_{h_n}$. Let $K \in  \mathcal{K}_{h_n}$, if $K \not \in \mathcal{K}^i_{h_n}$ for $i=0,2,4$ then $K \in \mathcal{K}^1_{h_n}$ and/or $K \in \mathcal{K}^3_{h_n}$, because $0 \le v_n \le 1$; hence, the union of the families $\mathcal{K}^i_{h_n}$ for $i=0,...,4$ is the whole $\mathcal{K}_{h_n}$. Moreover, if $K \in \mathcal{K}^i_{h_n}$ for $i=0,2,4$ then $K \not \in (\mathcal{K}^1_{h_n} \cup \mathcal{K}^3_{h_n})$, hence the sets $(\mathcal{K}^0_{h_n} \cup \mathcal{K}^2_{h_n} \cup \mathcal{K}^4_{h_n})$ and $(\mathcal{K}^1_{h_n} \cup \mathcal{K}^3_{h_n})$ are disjoint. It is clear, from the definition, that the families $\mathcal{K}^i_{h_n}$ are pairwise disjoint for $i=0,2,4$ because the corresponding sets $\Omega^i_{h_n}$ are pairwise disjoint. It remains to check that $\mathcal{K}^1_{h_n}$ and  $\mathcal{K}^3_{h_n}$ are disjoint, at least for $n \gg 1$.
Remember that $\| \nabla v_n \|_{L^\infty (\Omega)} \le C / \eps_n$, that $\mathrm{diam} (\tilde{K}) \le \tilde{C} h_n$ and that $h_n = o(\eps_n)$, then for $n \gg 1$ we have 
\begin{equation}   \label{e.stim}
	\begin{cases}  v_n  < 2 c_n + \bar{C} h_n / \eps_n < 1/3 \ \text{ in $\tilde K$}  &  \text{if $K \in \mathcal{K}^1_{h_n}$},  \\
		v_n \ge (1 - 2 c_n) - \bar{C} h_n / \eps_n  \ge 2/3  \ \text{ in $\tilde K$} &  \text{if $K \in \mathcal{K}^3_{h_n}$}. \end{cases}
\end{equation}
It follows that $\mathcal{K}^1_{h_n}$ and  $\mathcal{K}^3_{h_n}$ are disjoint for $n \gg 1$. 

Next, denote by $A^i_{h_n}$ the union of the elements $K \in \mathcal{K}^i_{h_n}$ and by $\tilde A^i_{h_n}$ the corresponding union of the extended elements (see Figure \ref{fig3}). Since the sets $\mathcal{K}^i_{h_n}$ provide a disjoint partition of $\mathcal{K}_{h_n}$ the corresponding sets $A^i_{h_n}$ give a disjoint partition of $\Omega$.
We claim that  for $n \gg 1$ the sets $A^i_{h_n}$ provide a disjoint partition of $\Omega$ (up to a set of measure zero, given by the union of the boundaries of the elements $K \in \mathcal{K}_{h_n}$). Moreover,
\begin{gather} 
	w_{h_n} = 0 \text{ in $A^0_{h_n}$} , \qquad c_n \le w_{h_n} \le 1 -c_n  \text{ in $A^2_{h_n}$} ,  \qquad w_{h_n} = 1 \text{ in $A^4_{h_n}$} , 
	\label{e.mn} \\
	-c_n \le w_{h_n} \le 1/3 + c_n  \text{ in $\tilde A^1_{h_n}$} , \qquad 2/3 - c_n \le w_{h_n} \le 1+ c_n  \text{ in $\tilde A^3_{h_n}$}  \label{e.ln}\\
	0 \le w_{h_n} \le 1/3 + c_n  \text{ in $\tilde A^1_{h_n} \setminus A^1_{h_n}$} , \qquad 2/3 - c_n \le w_{h_n} \le 1 \text{ in $\tilde A^3_{h_n} \setminus A^3_{h_n}$} . 
\label{e.li}
 \end{gather}
%
Let us check \eqref{e.mn}. By definition, if $K \in \mathcal{K}^0_{h_n}$  then $v_n =0$ on $\tilde K$, hence $w_{h_n} = 0$ in $K$ because the projection operator is locally an identity for constant functions (see for instance  \cite[Lemma 3.2]{BazilevsBeiraoCottrellHughesSangalli_M3AS06}). In the same way,  if $K \in \mathcal{K}^4_{h_n}$  then $v_n =1 $ on $\tilde K$, hence $w_{h_n} = 1$ in $K$. If $K \in \mathcal{K}^2_{h_n}$  then $ 2 c_n \le v_n \le 1 - 2 c_n$ in $K$ and then by \eqref{e.linfty} we have $ c_n \le w_{h_n} \le 1 - c_n$ in $K$. 
Let us check \eqref{e.ln}. If $K \in \mathcal{K}_{h_n}^1$ then, being $v_n \ge 0$, by \eqref{e.stim} we have $0 \le v_n  < 1/3$ in $\tilde{K}$ and thus by \eqref{e.linfty} we get $- c_n \le w_{h_n} \le 1/3 + c_n$ in $\tilde K$. Similarly for $K \in \mathcal{K}_{h_n}^3$. To get \eqref{e.li} from \eqref{e.ln} it is enough to note that $\tilde A^1_{h_n} \setminus A^1_{h_n}$ is contained in the union of the set $A^i_{h_n}$ for $i=0,2,3,4$ where $w_{h_n} \ge 0$. Similarly for $K \in \mathcal{K}_{h_n}^3$. 
Finally, note that \eqref{e.ln} implies that $\tilde{A}^1_{h_n}$ and $\tilde{A}^3_{h_n}$ are disjoint.

\separe

%
We are now ready to modify the function $w_{h_n}$ in the sets $A^i_{h_n}$ for $i=1,3$ (where the constraint $0 \le w_{h_n} \le 1$ may not be  satisfied).
Consider all the basis functions $\hat v_{h_n}$ 
whose support intersects an element $K \in \mathcal{K}^1_{h_n}$ and 
denote by $v^1_{h_n}$ their sum. 
%
By definition, basis functions $\hat v_{h_n}$ are non-negative, provide locally (on each element) a partition of unity and are supported in the extended elements $\tilde K$; hence 
\begin{gather*}
	0 \le v^1_{h_n}  \le 1  \text{ in $\Omega$} , \qquad v^1_{h_n} = 1 \text{ in $A^1_{h_n}$},  
\qquad \mathrm{supp} ( v^1_{h_n} ) \subset \tilde A^1_{h_n}, \\  \| D^{k} v^1_{h_n} \|_{L^\infty (\Omega)} \le C h_n^{-k} \ \text{ for $0 \le k \le 2$}.
\end{gather*}
The $L^\infty$-estimate for the derivatives follows from scaling and from the fact that the parametrization map $\mathbf{F}: [0,1]^3 \mapsto \Omega$ is 
a diffeomorphism of class $W^{2,\infty}$. Note that the support is contained in the enlarged set $\tilde{A}^1_{h_n}$.
Similarly we define $v^3_{h_n}$ and finally we can introduce the B-spline $v_{h_n}$, given by 
$$ 
     v_{h_n} = w_{h_n} + c_n v^1_{h_n} - c_n v^3_{h_n} .
$$
Since the supports of $v^1_{h_n}$ and $v^3_{h_n}$ are the disjoint sets $\tilde{A}^1_{h_n}$ and $\tilde{A}^3_{h_n}$, we can write $v_{h_n}$ as 
$$
	v_{h_n} = \begin{cases}
	w_{h_n} + c_n v^1_{h_n}	& \text{in } \tilde A^1_{h_n}	, \\
	w_{h_n} - c_n v^3_{h_n}  & \text{in } \tilde A^3_{h_n},	 \\
	w_{h_n}  & \text{otherwise}. \\
	\end{cases}
$$
 In the set $A^1_{h_n}$ we have $v^1_{h_n}  = 1$ and $-c_n \le w_{h_n} \le 1/3 +c_n $, hence $0 \le v_{h_n} \le 1$. In $ \tilde A^1_{h_n} \setminus  A^1_{h_n}$ we have $0 \le v^1_{h_n}  \le 1$ and $ 0 \le w_{h_n} \le 1-c_n$, hence $0 \le v_{h_n} \le 1$. 
%
%
%
We can argue in a similar way for $v^3_{h_n}$. We have checked that that $0 \le v_{h_n} \le 1$ in $\Omega$, for $n \gg 1$. 
%
%
Now, let us provide some error estimates. Writing 
$$
	v_{h_n} - w_{h_n} = \begin{cases}
	c_n v^1_{h_n}	& \text{in } \tilde A^1_{h_n}	, \\
	c_n v^3_{h_n}  & \text{in } \tilde A^3_{h_n},	 \\
	0 & \text{otherwise}. \\
	\end{cases}
$$
and using the $L^\infty$-estimates on $D^k v^1_{h_n}$ and $D^k v^3_{h_n}$we get 
$$  
\| D^k ( v_{h_n} - w_{h_n} ) \|_{L^\infty (\Omega)} \le C c_n h_n^{-k} \quad  \text{ for $0 \le k \le 2$} .$$ 

Let us check that the Lebesgue measure of $\mathrm{supp} (  v_{h_n} - w_{h_n} )$ is of order $\eps_n$. Clearly $\mathrm{supp} (  v_{h_n} - w_{h_n} ) \subset  ( \tilde A^1_{h_n} \cup \tilde A^3_{h_n})$. 
By Proposition \ref{p.Glsup} $ v_n = 1 $ in $\Omega \setminus (A_n \cup B_n)$ where $A_n \cup B_n = J_{2 \delta} \times (-\eps_n^\sharp , \eps_n^\sharp)$ where $J_{2\delta} = (1+2\delta) J$. 
Thus the sets $\Omega^1_n = \{ 0 < v_n < 2 c_n \} $ and $\Omega^3_n = \{ 1- 2 c_n < v_n < 1 \} $ are contained in $A_n \cup B_n$. It follows that $\tilde A^1_{h_n}$ and $\tilde A^3_{h_n}$ are contained in a set of the form   
$$ (1 + 2\delta + \tilde{C} h_n) J \times ( - \eps_n^\sharp - \tilde{C} h_n , \eps_n^\sharp + \tilde{C} h_n) . $$
Since $h_n = o (\eps_n)$ we have the required estimate on the measure of the support. Then, using $ c_n = C (h_n/\eps_n)^3$ and the $L^\infty$-estimates above we get  
\begin{align}
	& \int_\Omega | v_{h_n} - w_{h_n} |^2 \, dx \le C c_n^2 \eps_n \le C' h_n^6 \eps_n^{-5} ,  \label{e.L2ee} \\
%
	& \int_\Omega | \nabla v_{h_n} - \nabla w_{h_n}   |^2 \le C c_n^2 \eps_n h_n^{-2} \le C' h_n^4 \eps_n^{-5} ,  \label{e.H1ee} \\
	& \int_\Omega | D^2 v_{h_n} - D^2  w_{h_n}  |^2  \le C c_n^2 \eps_n h_n^{-4} \le C'  h_n^2 \eps_n^{-5} . \label{e.H2ee} \phantom{-------}
\end{align}
Before proceeding, let us provide also some global error estimates. We know (see for instance  \cite[Lemma 3.2]{BazilevsBeiraoCottrellHughesSangalli_M3AS06}) that $\Pi_{\V_{h_n}} 1 = 1$ and then
$v_n - w_{h_n}  =  (v_n -1)  - \Pi_{\hspace{0pt}\V_{h_n}} (v_n -1) $.  Hence, using \eqref{e.UVhOm} for $(v_n -1)$ and $l=3$ together with Corollary  \ref{c.corlsup} we get, for $k=0,..,2$,
\begin{equation} \label{e.19}
	| v_n  - w_{h_n}|^2_{H^k (\Omega)} \le C \,  h_n^{6-2k} \| v_n -1  \|^2_{H^3 (\Omega)}  \le C' h_n^{6-2k} \eps_n^{-5} .
\end{equation}
Note that these estimate for $k=0,...,2$ sre of the same order of \eqref{e.L2ee}-\eqref{e.H2ee}.

\separe

{\bf Step 2.} Now, let us prove \eqref{e.18h}. In the sequel we will make frequent use of the following Young's inequality $(a+ b)^2 \le (1+\delta^{-1}) a^2 + ( 1 + \delta) b^2$ for $\delta >0$. Let $C'_\delta = (1+\delta)(1+\delta^{-1})$. Using twice Young's inequality, the error estimates \eqref{e.L2ee} and \eqref{e.19} we get 
\begin{align*}
	\int_\Omega \eps_n^{-1} | v_{h_n} - 1 |^2 \, 
	& \le 
	(1+\delta^{-1}) \int_\Omega \eps_n^{-1} | v_{h_n} - w_{h_n} |^2 \, dx  +  ( 1 + \delta)  \int_\Omega \eps_n^{-1} | w_{h_n} - 1 |^2 \, dx  \\
	& \le
	(1+\delta^{-1}) \int_\Omega \eps_n^{-1} | v_{h_n} - w_{h_n} |^2 \, dx + C'_\delta   \int_\Omega \eps_n^{-1} | w_{h_n} - v_n |^2 \, dx  \, + \\
	& \quad + ( 1 + \delta)^2   \int_\Omega \eps_n^{-1} | v_{n} - 1 |^2 \, dx  \\ 
	& \le 
	C_\delta \, (h_n / \eps_n)^6 + ( 1 + \delta)^2   \int_\Omega \eps_n^{-1} | v_{n} - 1 |^2 \, dx .
\end{align*}
Similarly, 
\begin{align*} 
	\int_\Omega \eps_n  | \nabla v_{h_n} |^2 \, dx 
	& \le 
	(1+ \delta^{-1} )  \int_\Omega \eps_n | \nabla v_{h_n} -  \nabla w_{h_n}  |^2 \, dx +  (1+\delta) \int_\Omega  \eps_n | \nabla w_{h_n} |^2 \, dx  \\
	& \le 
	(1+ \delta^{-1} )  \int_\Omega \eps_n | \nabla v_{h_n} -  \nabla w_{h_n}  |^2 \, dx + C'_\delta \int_\Omega  \eps_n | \nabla w_{h_n} - \nabla v_n|^2 	\, dx \\ &  \quad + ( 1 + \delta)^2 \int_\Omega  \eps_n | \nabla v_n |^2 \, dx  \\
	& \le  C_\delta \, (h_n / \eps_n)^4   + ( 1 + \delta)^2 \int_\Omega  \eps_n | \nabla v_n |^2 \, dx .
\end{align*}
Finally, 
\begin{align*} 
	\int_\Omega \eps^3_n  | \Delta v_{h_n} |^2 \, dx 
	& \le (1+ \delta^{-1} )  \int_\Omega \eps^3_n | \Delta v_{h_n} -  \Delta w_{h_n}  |^2 \, dx +  (1+\delta) \int_\Omega  \eps^3_n | \Delta w_{h_n} |^2 \, dx  \\	& \le 
	(1+ \delta^{-1} )  \int_\Omega \eps^3_n | \Delta v_{h_n} -  \Delta w_{h_n}  |^2 \, dx + C'_\delta \int_\Omega  \eps^3_n | \Delta w_{h_n} - \Delta v_n|^2 	\, dx \\ &  \quad + ( 1 + \delta)^2 \int_\Omega  \eps^3_n | \Delta v_n |^2 \, dx  \\
	& \le  C_\delta \, (h_n / \eps_n)^2  + ( 1 + \delta)^2 \int_\Omega  \eps_n | \Delta v_n |^2 \, dx .
\end{align*}
In conclusion,
\begin{align*}
	 \int_{\Omega}  \eps_n^{-1} |v_{h_n} -1|^2 + 2 \eps_n | \nabla v_{h_n} |^2 & + \eps_n^3 | \Delta v_{h_n} |^2  \, dx  \le \\ & \le  (1 +\delta)^2
	\int_{\Omega}  \eps_n^{-1} |v_{n} -1|^2 + 2 \eps_n | \nabla v_{n} |^2 + \eps_n^3 | \Delta v_{n} |^2  \, dx + o (1)
\end{align*}
and thus, by Proposition \ref{p.Glsup} 
$$
	 \limsup_{n \to +\infty} \int_{\Omega}  \eps_n^{-1} |v_{h_n} -1|^2 + 2 \eps_n | \nabla v_{h_n} |^2  + \eps_n^3 | \Delta v_{h_n} |^2  \, dx  \le 
	4 | J | + C \delta ,
$$
for a suitable $C>0$, depending only on $J$.

\separe

{\bf Step 3.} Let $u_{h_n} = \Pi_{\U_{h_n}} u_n $. 
Since the interpolation estimates in $\U_{h_n}$ and $\V_{h_n}$ are non-local it is necessary to introduce a further set, ``between'' $E_n = J_{\delta/2} \times (-\eps_n^\flat /2, \eps_n^\flat /2)$ and $C_n = J_\delta \times (-\eps_n^\flat, \eps_n^\flat)$    (see Corollary \ref{c.corlsup} and \ref{c.corlsupbis}): for $ \tfrac12 < a < 1$ let $J_n = J_{a \delta } \times (- a \eps_n^\flat, a \eps_n^\flat)$. 

%
%
%
%
%
%
Since $W$ is quadratic, by Young's inequality we can write 
$$
	W (\strain (u_{h_n}) ) = W \big( \strain (u) - \strain (u - u_{h_n}) \big) \le (1+ \delta) W ( \strain (u) ) + C_\delta | D u - D u_{h_n} |^2 .
$$
Being $v_{h_n} \le 1$ we get 
\begin{align*}
\int_{\Omega \setminus J_n} ( v_{h_n}^2 +\eta_n) W ( \strain(u_{h_n}) )  \le & \ (1+ \eta_n)   \int_{\Omega \setminus J_n}  (1+ \delta)   W ( \strain (u) ) + C_\delta  | D u - D u_{h_n} |^2\, dx .
\end{align*}
Clearly, for the first term we have
$$
  \limsup_{n \to +\infty} \ (1+ \eta_n)  \int_{\Omega \setminus J_n}  W ( \strain (u) )  \, dx \le   
  \int_{\Omega \setminus J}  W ( \strain (u) ) \, dx .
$$
From Corollary \ref{c.corlsupbis}  we know that $u_n = u$ in $\Omega \setminus E_n$ where $E_n = J_{\delta/2} \times (-\eps_n^\flat/2, \eps_n^\flat/2)$. Moreover $u \in W^{2,\infty} ( \Omega \setminus J)$. Since $h_n = o( \eps_n)$ we can employ the (non-local) interpolation error estimate \eqref{e.err-Uhloc}, i.e. 
$$
	|  u_n - \Pi_{\hspace{1.2pt}\U_{h_n}} u_n |_{H^k(K,\,\R^3)} \le C h_n^{2-k} \| u_n \|_{H^2 (\tilde{K}, \,\R^3)} 
$$
to obtain
$$
\int_{\Omega \setminus J_n}   | D u_n - D u_{h_n} |^2\, dx 
\le C h^2_n \to 0. 
$$
%
%
%
%
%
From Corollary \ref{c.corlsup}  we know that $v_n = 0$ in $C_n = J_\delta \times ( -\eps_n^\flat, \eps_n^\flat)$, hence $v_{h_n} = 0$ in $J_n$ and we can write 
\begin{align*}
\int_{J_n} ( v_{h_n}^2 +\eta_n) W ( \strain(u_{h_n}) ) \, dx  & \le C \eta_n \int_{J_n} | Du_{h_n} |^2 \, dx 
\\
& \le C' \eta_n  \int_{J_n} | Du_n |^2 \, dx +  C' \eta_n \int_{J_n} | Du_n -  Du_{h_n}|^2 \, dx .
\end{align*}
By Corollary \ref{c.corlsupbis} we know that $\| D u_n \|_{L^\infty ( \Omega )} \le C  \eps_n^{-1}$. Being $\eta_n = o(\eps_n)$ we get 
$$
	\eta_n  \int_{J_n} | Du_n |^2 \, dx \le C \eta_n | J _\delta | \eps_n^{-1} \to  0 .
$$
Moreover, the error estimate \eqref{e.err-Uhloc} for $k=l=1$ provides
$$
	|  u_n - \Pi_{\hspace{1.2pt}\U_{h_n}} u_n |^2_{H^1(K,\,\R^3)} \le C \| u \|^2_{H^1 (\tilde{K}, \,\R^3)} \le C h_n^3  \eps_n^{-2} .
$$
Hence, for $\eta_n = o(\eps_n)$, 
$$
	 \eta_n \int_{J_n} | Du_n -  Du_{h_n}|^2 \, dx \le C \eta_n \, ( | J_n | / h_n^3 ) \, h_n^3 \eps_n^{-2}
	\le C' \eta_n  \eps_n^\flat \eps_n^{-2} \le C'' \eta_n \eps_n^{-1}  \to 0 .
$$
The proof is concluded. \qed

\begin{remark} \label{r.remo} Note that having the estimate $-c_n \le w_{h_n} \le 1 + c_n$ for $c_n = (h_n / \eps_n)^3$ it is not possible to employ the linear transform $v_{h_n} = \ell ( w_{h_n} )$ where $\ell(t) = (t + c_n) /(1 + 2 c_n)$. Indeed, $0 \le v_n \le 1$ but 
$$
     \int_{ \{ w_{h_n} \ge 1 \} } \eps_n^{-1} | v_{h_n}  -1 |^2 \, dx \ge C \eps^{-1} c_n^2 
$$ 
and $\eps_n^{-1} c_n^2$ is even not bounded under the assumption $c_n = o (1)$, i.e.~$h_n = o (\eps_n)$. Possibly this simply construction could work under more restrictive assumptions on the ratio between the mesh size $h_n$ and the internal length $\eps_n$. 
\end{remark} 



\appendix 

\section{{\boldmath$GSBD$ spaces} \label{A}}

We provide just the definition and the main properties of vector fields in $GSBD (\Omega)$ and $GSBD^2(\Omega)$ for $\Omega$ an open subset of $\R^3$. For a general and detailed work the reader should refer to \cite{DalMaso_JEMS13}. 

\medskip
For $\xi \in \mathbb{S}^2 = \{ \xi \in \R^3 : | \xi | =1 \}$ let $\xi^\perp = \{ y \in \R^3 : y \cdot \xi =0 \}$. For $B \subset \Omega$ and $y \in \xi^\perp$ let $B^\xi_y = \{ s \in \mathbb{R} : y + s \xi  \in  B\}$ denote the ``slice" of $B$. If $u: \Omega \to \R^3$ we consider its projected $\xi$-slice in $B$, i.e., the function $u^\xi_y : B^\xi_y \to \mathbb{R}$ given by $ u^\xi_y ( s) = u ( y + s \xi ) \cdot \xi$.
Note that $u^\xi_y$ is scalar valued.


\begin{definition} \label{d.GSBDdef} A measurable function $u : \Omega \to \R^3$ belongs to $GSBD(\Omega, \R^3)$ if for every $\xi \in \mathbb{S}^2$ and a.e.~$y \in \xi^\perp$ the slices $u^\xi_y$ belong to $SBV_{loc} ( \Omega^\xi_y)$ and if there exists a finite Radon measure $\mu$ such that for every Borel set $B \subset \Omega$ we have 
$$
	\int_{\xi^\perp} | D  u^\xi_y |  \big( B^\xi_y \setminus J^{(1)} ( u^\xi_y ) \big) +  \# \big( B^\xi_y \cap J^{(1)}  ( u^\xi_y ) \big) \, dy \le  \mu (B)  \quad  \text{for every $\xi \in \mathbb{S}^2$ and a.e.~$y \in \xi^\perp$.}  
$$
Here $D u^\xi_y \in \mathcal{M}_{loc}(\Omega^\xi_y)$   is the distributional derivative of $u^\xi_y$ while $ J^{(1)}  ( u^\xi_y ) = \{  s \in \Omega^\xi_y : | \llbracket u^\xi_y (s) \rrbracket | \ge 1 \}$.
\end{definition}

\begin{theorem} \label{t.GSBDslic} Let $u \in GSBD (\Omega)$ and $\xi \in \mathbb{S}^2$. For a.e.~$y \in \xi^\perp$ we have  $(J^\xi (u) )_y^\xi = J (u^\xi_y)$ where 
$$
	J^\xi (u) = \{ x \in J (u)  :  (u^+ (x)  - u^-(x) ) \cdot \xi \neq 0  \}  .
$$
Moreover, for a.e.~$\xi \in \mathbb{S}^2$ we have $\H^2 ( J^\xi (u) \setminus J (u) ) = 0$ and 
$$
	\int_{\xi^\perp} \# ( J (u^\xi_y) ) \, dy = \int_{J(u)} | \xi \cdot \nu | \, d \mathcal{H}^2 .
$$
\end{theorem}



\begin{definition} \label{d.GSBD2def} A measurable function $u : \Omega \to \R^3$ belongs to $GSBD^2(\Omega, \R^3)$ if $u \in GSBD(\Omega, \R^3)$, $\strain(u) \in L^2 ( \Omega , \mathbb{R}^{ 3 \times 3})$ and $\mathcal{H}^2 (J(u)) < +\infty$. 
\end{definition}

Combining \cite{Iurlano_CVPDE12} and \cite{CortesaniToader_NA99} yields the following approximation result.

\begin{theorem} \label{t.GSBDdens} Let $u \in GSBD^2 (\Omega) \cap L^2 (\Omega, \R^3)$. Then  there exists a sequence 
$u_k \in SBV^2 (\Omega , \R^3)$ such that  
$u_k \to u$ in $L^2 (\Omega, \R^3)$,
$\strain (u_k) \to \strain (u)$ in $L^2 (\Omega, \R^{3 \times 3})$
and $\mathcal{H}^2  ( J (u_k) ) \to  \mathcal{H}^2 ( J(u))$.
Further, $u_k$ can be chosen in such a way that 
$J (u_k) \subset \Omega$ is the finite union of closed, disjoint simplexes
and $u_k \in W^{m,\infty} ( \Omega \setminus J (u_k) , \R^3)$ (for $m$ arbitrarily large).
\end{theorem}

{\bf Aknowledgement.} The author wishes thank G.~Sangalli and A.~Bressan for helpful discussions on isogeometric B-splines.
Financial support has been provided by GNAMPA-INdAM project ``Analisi multiscala di sistemi complessi
con metodi variazionali"  and by ERC Advanced Grant ``QuaDynEvoPro" \#290888.

\bibliographystyle{plain}
\bibliography{gamma} 

\begin{thebibliography}{10}

\bibitem{AlmiNegri_19}
S.~Almi and M.~Negri.
\newblock Analysis of staggered evolutions for nonlinear energies in phase
  field fracture.
\newblock {\em arXiv:1904.01895}.

\bibitem{AmbatiGerasimovDeLorenzi_CM15}
M.~Ambati, T.~Gerasimov, and L.~{De Lorenzis}.
\newblock {A review on phase-field models of brittle fracture and a new fast
  hybrid formulation}.
\newblock {\em Comp. Mech.}, 55(2):383--405, 2015.

\bibitem{AmbrosioTortorelli_CPAM90}
L.~Ambrosio and V.M. Tortorelli.
\newblock {Approximation of functionals depending on jumps by elliptic
  functionals via {$\Gamma$}-convergence}.
\newblock {\em Comm. Pure Appl. Math.}, 43(8):999--1036, 1990.

\bibitem{ArtinaFornasierMichelettiPerotto_SIAMJSC15}
M.~Artina, M.~Fornasier, S.~Micheletti, and S.~Perotto.
\newblock Anisotropic mesh adaptation for crack detection in brittle materials.
\newblock {\em SIAM J. Sci. Comput.}, 37(4):B633--B659, 2015.

\bibitem{Bach_ESAIMCOCV18}
A.~Bach.
\newblock Anisotropic free-discontinuity functionals as the {$\Gamma$}-limit of
  second-order elliptic functionals.
\newblock {\em ESAIM Control Optim. Calc. Var.}, 24(3):1107--1140, 2018.

\bibitem{BazilevsBeiraoCottrellHughesSangalli_M3AS06}
Y.~Bazilevs, L.~Beir{\~a}o~da Veiga, J.A. Cottrell, T.J.R. Hughes, and
  G.~Sangalli.
\newblock Isogeometric analysis: approximation, stability and error estimates
  for {$h$}-refined meshes.
\newblock {\em Math. Models Methods Appl. Sci.}, 16(7):1031--1090, 2006.

\bibitem{BeiraoBuffaSangalliVazquez_16}
L.~Beir\~ao~da Veiga, A.~Buffa, G.~Sangalli, and R.~V\'azquez.
\newblock An introduction to the numerical analysis of isogeometric methods.
\newblock In {\em Numerical simulation in physics and engineering}, volume~9 of
  {\em SEMA SIMAI Springer Ser.}, pages 3--69. Springer, 2016.

\bibitem{BellettiniCoscia_NFAO94}
G.~Bellettini and A.~Coscia.
\newblock Discrete approximation of a free discontinuity problem.
\newblock {\em Numer. Funct. Anal. Optim.}, 15(3-4):201--224, 1994.

\bibitem{BordenHughesLandisVerhoosel_CMAME14}
M.J. Borden, T.J.R. Hughes, C.M. Landis, and C.V. Verhoosel.
\newblock A higher-order phase-field model for brittle fracture: Formulation
  and analysis within the isogeometric analysis framework.
\newblock {\em Comput. Methods Appl. Mech. Engrg.}, 273:100 -- 118, 2014.

\bibitem{BourdinChambolle_NM00}
B.~Bourdin and A.~Chambolle.
\newblock Implementation of an adaptive finite-element approximation of the
  {M}umford-{S}hah functional.
\newblock {\em Numer. Math.}, 85(4):609--646, 2000.

\bibitem{BourdFrancMar00}
B.~Bourdin, G.~A. Francfort, and J.-J. Marigo.
\newblock {Numerical experiments in revisited brittle fracture}.
\newblock {\em J. Mech. Phys. Solids}, 48(4):797--826, 2000.

\bibitem{BourdFrancMar08}
B.~Bourdin, G.A. Francfort, and J.-J. Marigo.
\newblock The variational approach to fracture.
\newblock {\em J. Elasticity}, 91:5--148, 2008.

\bibitem{Braides98}
A.~Braides.
\newblock {\em {Approximation of free-discontinuity problems}}.
\newblock Springer-Verlag, Berlin, 1998.

\bibitem{BurgerEspositoZeppieri_MMS15}
M.~Burger, T.~Esposito, and C.I. Zeppieri.
\newblock Second-order edge-penalization in the {A}mbrosio-{T}ortorelli
  functional.
\newblock {\em Multiscale Model. Simul.}, 13(4):1354--1389, 2015.

\bibitem{BurkOrtnSuel10}
S.~Burke, C.~Ortner, and E.~S{\"u}li.
\newblock {An adaptive finite element approximation of a variational model of
  brittle fracture}.
\newblock {\em SIAM J. Numer. Anal.}, 48(3):980--1012, 2010.

\bibitem{Chambolle_JMPA04}
A.~Chambolle.
\newblock An approximation result for special functions with bounded
  deformation.
\newblock {\em J. Math. Pures Appl. (9)}, 83(7):929--954, 2004.

\bibitem{CCF}
A.~Chambolle, S.~Conti, and G.A. Francfort.
\newblock Approximation of a brittle fracture energy with a constraint of
  non-interpenetration.
\newblock {\em Arch. Ration. Mech. Anal.}, 228(3):867--889, 2018.

\bibitem{ChambolleCrismale_ARMA19}
A.~Chambolle and V.~Crismale.
\newblock A density result in {$GSBD^p$} with applications to the approximation
  of brittle fracture energies.
\newblock {\em Arch. Ration. Mech. Anal.}, 232(3):1329--1378, 2019.

\bibitem{Ciarl78}
P.G. Ciarlet.
\newblock {\em {The finite element method for elliptic problems}}.
\newblock {Studies in Mathematics and its Applications, Vol. 4}. North-Holland
  Publishing Co., Amsterdam, 1978.

\bibitem{CortesaniToader_NA99}
G.~Cortesani and R.~Toader.
\newblock {A density result in {SBV} with respect to non-isotropic energies}.
\newblock {\em Nonlinear Anal.}, 38(5):585--604, 1999.

\bibitem{DalMaso93}
G.~{Dal Maso}.
\newblock {\em {An introduction to {$\Gamma$}-convergence}}.
\newblock Birkh{\"a}user, Boston, 1993.

\bibitem{DalMaso_JEMS13}
G.~{Dal Maso}.
\newblock {Generalised functions of bounded deformation}.
\newblock {\em J. Eur. Math. Soc.}, 15(5):1943--1997, 2013.

\bibitem{DalMasoIurlano_CPAA13}
G.~{Dal Maso} and F.~Iurlano.
\newblock Fracture models as {$\Gamma$}-limits of damage models.
\newblock {\em Commun. Pure Appl. Anal.}, 12(4):1657--1686, 2013.

\bibitem{FonsecaMantegazza_SJMA00}
I.~Fonseca and C.~Mantegazza.
\newblock Second order singular perturbation models for phase transitions.
\newblock {\em SIAM J. Math. Anal.}, 31(5):1121--1143, 2000.

\bibitem{Iurlano_CVPDE12}
F.~Iurlano.
\newblock A density result for {$GSBD$} and its application to the
  approximation of brittle fracture energies.
\newblock {\em Calc. Var. Partial Differential Equation}, 51:315--342, 2014.

\bibitem{KiendlAmbatiDeLorenzisGomezReali_CMAME16}
J.~Kiendl, M.~Ambati, L.~De Lorenzis, H.~Gomez, and A.~Reali.
\newblock Phase-field description of brittle fracture in plates and shells.
\newblock {\em Comput. Methods Appl. Mech. Engrg.}, 312:374 -- 394, 2016.

\bibitem{KneesNegri_15}
D.~Knees and M.~Negri.
\newblock Convergence of alternate minimization schemes for phase field
  fracture and damage.
\newblock {\em Math. Models Methods Appl. Sci.}, 27(9):1743--1794, 2017.

\bibitem{KuhnMueller_EFM10}
C.~Kuhn and R.~M\"uller.
\newblock A continuum phase field model for fracture.
\newblock {\em Engineering Fracture Mechanics}, 77(18):3625 -- 3634, 2010.

\bibitem{LarsenOrtnerSuli_M3AS10}
C.J. Larsen, C.~Ortner, and E.~S{\"u}li.
\newblock {Existence of solutions to a regularized model of dynamic fracture}.
\newblock {\em Math. Models Methods Appl. Sci.}, 20(7):1021--1048, 2010.

\bibitem{LiPecoMillanAriasArroyo_IJNME15}
B.~Li, C.~Peco, D.~Mill\'an, I.~Arias, and M.~Arroyo.
\newblock Phase-field modeling and simulation of fracture in brittle materials
  with strongly anisotropic surface energy.
\newblock {\em Int. J. Numer. Methods Eng.}, 102(3-4):711--727, 2015.

\bibitem{March_VC92}
R.~March.
\newblock Visual reconstruction with discontinuities using variational methods.
\newblock {\em Vision Computing}, 10:30--38, 1992.

\bibitem{MieheWelschingerHofacker10}
C.~Miehe, F.~Welschinger, and M.~Hofacker.
\newblock {Thermodynamically consistent phase-field models of fracture:
  variational principles and multi-field {FE} implementations}.
\newblock {\em Internat. J. Numer. Methods Engrg.}, 83(10):1273--1311, 2010.

\bibitem{MikelicWheelerWick_N15}
A.~Mikeli{\'c}, M.~F. Wheeler, and T.~Wick.
\newblock A quasi-static phase-field approach to pressurized fractures.
\newblock {\em Nonlinearity}, 28(5):1371--1399, 2015.

\bibitem{N_LNACM16}
M.~Negri.
\newblock Quasi-static evolutions in brittle fracture generated by gradient
  flows: sharp crack and phase-field approaches.
\newblock In {\em Innovative Numerical Approaches for Multi-Physics and
  Multi-Scale Problems}, volume~81 of {\em Lecture Notes in Applied and
  Computational Mechanics}, pages 197--216. Springer, 2016.

\bibitem{PaulZimmermannMandadapuHughesLandisSauer_19}
K.~Paul, C.~Zimmermann, K.K. Mandadapu, T.J.R. Hughes, C.M. Landis, and R.A.
  Sauer.
\newblock An adaptive space-time phase field formulation for dynamic fracture
  of brittle shells based on {LR} {NURBS}.
\newblock {\em arXiv:1906.10679}.

\bibitem{Schumaker_07}
L.L. Schumaker.
\newblock {\em Spline functions: basic theory}.
\newblock Cambridge Mathematical Library. Cambridge University Press,
  Cambridge, third edition, 2007.

\bibitem{SicsicMarigoMaurini_JMPS14}
P.~Sicsic, J.-J. Marigo, and C.~Maurini.
\newblock Initiation of a periodic array of cracks in the thermal shock
  problem: A gradient damage modeling.
\newblock {\em J. Mech. Phys. Solids}, 63:256--284, 2014.

\end{thebibliography}

\end{document}